\pdfoutput=1
\nonstopmode
\documentclass[]{amsart}
\usepackage{amssymb}
\usepackage{amscd}
\usepackage{graphicx}
\usepackage{epsfig}
\usepackage[all]{xy}
\newdir{ >}{{}*!/-10pt/@{>}}
\newdir^{ (}{{}*!/-5pt/@^{(}}
\newdir_{ (}{{}*!/-5pt/@_{(}}
\def\cgaps#1{}
\def\Cgaps#1{}
\allowdisplaybreaks

\def\undersetbrace#1\to#2{\underbrace{#2}_{#1}}
\def\oversetbrace#1\to#2{\overbrace{#2}^{#1}}
\def\AMSunderset#1\to#2{\underset{#1}{#2}}
\def\AMSoverset#1\to#2{\overset{#1}{#2}}

\def\East#1#2{\xrightarrow[\;\;#2\;\;]{\;\;#1\;\;}}

\newcommand{\nmb}[2]{\ifx!#1{\ref{nmb:#2}}%
\else\if.#1{\label{nmb:#2}}%
\else\if0#1{\label{nmb:#2}}%
\else{{#2}}%
\fi\fi\fi}
\swapnumbers

\newtheorem*{proposition*}{Proposition}
\newtheorem{theorem}[subsection]{Theorem}
\newtheorem*{theorem*}{Theorem}

\newtheorem*{lemma*}{Lemma}

\newtheorem*{corollary*}{Corollary}
\newenvironment{demo}[1]{\par\smallskip\noindent{\bf #1.}}{\par\smallskip}

\def\o{\circ\,}
\def\X{\mathfrak X}
\def\al{\alpha}
\def\be{\beta}
\def\ga{\gamma}
\def\de{\delta}

\def\et{\eta}

\def\ka{\kappa}
\def\la{\lambda}
\def\rh{\rho}

\def\ph{\varphi}

\def\om{\omega}
\def\Ga{\Gamma}
\def\De{\Delta}

\def\La{\Lambda}

\def\Ph{\Phi}

\def\Om{\Omega}
\def\i{^{-1}}
\def\x{\times}
\def\p{\partial}
\let\on=\operatorname

\def\L{\mathcal L}
\def\AMSonly#1{}
\def\s{^{\sharp}}
\def\sh{\sharp}

\def\ad{\on{ad}}
\def\R{\mathbb R}

\def\<{\big\langle}
\def\>{\big\rangle \:}
\def\Fl{\operatorname{Fl}}
\def\h{^\text{hor}}

\begin{document}

\title[Sobolev curvature]
{Sobolev Metrics on Diffeomorphism Groups and the Derived Geometry of Spaces of Submanifolds}
\author{Mario Micheli, Peter W. Michor and David Mumford}
\address{Mario Micheli:
MAP5, Universit\'e Paris Descartes 
45, rue des Saints P\`eres, 7\`eme \'etage 
75270 Paris Cedex 06, France
}
\email{mariomicheli@gmail.com}
\address{
Peter W. Michor:
Fakult\"at f\"ur Mathematik, Universit\"at Wien,
Nordbergstrasse 15, A-1090 Wien, Austria}
\email{Peter.Michor@univie.ac.at}
\address{
David Mumford:
Division of Applied Mathematics, Brown University,
Box F, Providence, RI 02912, USA}
\email{David\_{}Mumford@brown.edu}
\thanks{
MM was supported by ONR grant N00014-09-1-0256, 
PWM was supported by FWF-project 21030, 
DM was supported by 
NSF grant DMS-0704213,
and all authors where supported by 
NSF grant DMS-0456253.}
\dedicatory{Dedicated to I.R.~Shafarevich on the occasion of his 90th birthday} 
\date{\today}
\begin{abstract}
Given a finite dimensional manifold $N$, the group $\on{Diff}_{\mathcal S}(N)$ of diffeomorphism of $N$ which 
fall suitably rapidly to the identity, acts 
on the manifold $B(M,N)$ of submanifolds on $N$ of diffeomorphism type $M$ where $M$ is a compact 
manifold with $\dim M<\dim N$. For a right invariant weak Riemannian metric on $\on{Diff}_{\mathcal S}(N)$ induced 
by a quite general operator $L:\X_{\mathcal S}(N)\to \Ga(T^*N\otimes\on{vol}(N))$, we consider the induced 
weak Riemannian metric on $B(M,N)$ and we compute its geodesics and sectional curvature.
For that we derive a covariant formula for curvature in finite and infinite dimensions, we show how 
it makes O'Neill's formula very transparent, and we use it finally to compute sectional curvature 
on $B(M,N)$.
\end{abstract}

\maketitle

\section{Introduction}
It was 46 years ago that Arnold discovered an amazing link between Euler's equation for 
incompressible non-viscous fluid flow and geodesics in the group of volume preserving 
diffeomorphisms $\on{SDiff}(\R^n)$ under the $L^2$-metric \cite{Arnold66}. One goal in this paper 
is to extend his ideas to a large class of Riemannian metrics on the group of all diffeomorphisms 
$\on{Diff}_{\mathcal S}(N)$ falling suitably to the identity, of any finite dimensional manifold $N$. 
The resulting geodesic equations are 
integro-differential equations for fluid-like flows on $N$ determined by an initial velocity field. 
In previous papers \cite{MM07, MGeomEvol, MMM1}, we have looked at the special case where $N=\R^n$ 
and the metric is a sum of Sobolev norms on each component of the tangent vector but here we 
develop the formalism to work in a very general setting.

The extra regularity given by using higher order norms means that these metrics on the group of 
diffeomorphisms can induce a metric on many quotient spaces of the diffeomorphism group modulo a 
subgroup. This paper focuses on the space of submanifolds of $N$ diffeomorphic to some $M$ which we 
denote by $B(M,N)$. $\on{Diff}_{\mathcal S}(N)$ acts on $B(M,N)$ with open orbits, one for each isotopy type of 
embedding of $M$ in $N$. The spaces $B$ may be called the {\it Chow manifolds} of $N$ by analogy 
with the Chow varieties of algebraic geometry, or {\it non-linear Grassmannians} because of their 
analogy with the Grassmannian of linear subspaces of a projective space. The key point is that the 
metrics we study will descend to the spaces $B(M,N)$ so that the map $\on{Diff}_{\mathcal S}(N) 
\rightarrow B(M,N)$ (given by the group action on a base point) is a Riemannian submersion. 
Geodesics from one submanifold to another may be thought of as deformations of one into the other 
realized by a flow on $N$ of minimal energy.

In the special case where $M$ is a finite set of points, $B(M,N)$ is called the space of 
{\it landmark point sets} in $N$. This has been used extensively by statisticians for example and 
is the subject of our previous paper \cite{MMM1}. The case $B(S^1,\R^2)$ is the space of all simple 
closed plane curves and has been studied in many metrics, see~\cite{kushnarev:1,MM98, SM06} for example. This 
and the case $B(S^2,\R^3)$ of spheres in 3-space have had many applications to medical imaging, 
constructing optimal warps of various body parts or sections of body parts from one medical scan to 
another \cite{MilEtAl,zhang_s}.

The high point of Arnold's analysis was his determination of the sectional curvatures in the group 
of volume preserving diffeomorphisms. This has had considerable impact on the analysis of the 
stability and instability of incompressible fluid flow. A similar formula for sectional curvature 
of $B(M,N)$ may be expected to shed light on how stable or unstable geodesics are in this space, 
e.g. whether they are unique and effective for medical applications.

Computing this curvature required a new formula. In general, the induced inner product on the 
cotangent space of a submersive quotient is much more amenable to calculations than the inner 
product on the tangent space. The first author found a new formula for the curvature tensor of a 
Riemannian manifold which uses only derivatives of the former, the dual metric tensor. This result, 
{\it Mario's formula}, is proven in section 2. In this section we also define a new class of 
infinite dimensional Riemannian manifolds, {\it robust Riemannian manifolds} to which Mario's 
formula and our analysis of submersive quotients applies. We also obtain a transparent new proof of 
O'Neill's curvature fomula. This class of manifolds builds on the 
theory of {\it convenient} infinite dimensional manifolds, see \cite{KM97}. 
To facilitate readability this theory is summarized in an Appendix. 

In section 3, we describe diffeomorphism groups of a finite dimensional manifold $N$ 
consisting of diffeomorphisms which decrease suitably rapidly to the identity on $N$ if we move to 
infinity on $N$; only these admit charts and are a regular Lie groups. 
We shall denote by $\on{Diff}_{\mathcal S}(N)$ any of these 
groups in order to simplify notation, and by $\X_{\mathcal S}(N)$ the corresponding Lie algebra of suitably 
decreasing vector fields on $N$.
We introduce 
a very general class of Riemannian metrics given by a positive definite self-adjoint differential 
operator $L$ from the space of smooth vector fields on $N$ to the space of measure-valued 1-forms. 
This defines an inner product on vector fields $X,Y$ by:
$$ \langle X, Y \rangle_L = \int_N (LX, Y).$$
Note that $LX$ paired with $Y$ gives a measure on $N$ hence can be integrated without assuming $N$ 
carries any further structure. Under suitable assumptions, the inverse of $L$ is given by a kernel 
$K(x,y)$ on $N \times N$ with values in $p_1^*TN \otimes p_2^*TN$. We then describe the geodesic 
equation in $\on{Diff}_{\mathcal S}(N)$ for these metrics. It is especially simple written in terms 
of the {\it momentum}. If $\ph(t)\in \on{Diff}_{\mathcal S}(N)$ is the geodesic, then 
$X(t) = \partial_t(\ph)\circ \ph\i$ is a time varying vector field on $N$ and its momentum is 
simply $LX(t)$. 

In section 4 we introduce the induced metrics on $B(M,N)$. We give the geodesic equation for these metrics also using momentum. One of the keys to working in this space is to define a convenient set of vector fields and forms on $B$ in terms of auxiliary forms and vector fields on $N$. In this way, differential geometry on $B$ can be reduced to calculations on $N$. Lie derivatives on $N$ are especially useful here. 

In the final section 5, we compute the sectional curvatures of $B(M,N)$. Like Arnold's formula, we get a formula with several terms each of which seems to play a different role. The first involves the second derivatives of $K$ and the others are expressed in terms which we call {\it force} and {\it stress}. Force is the bilinear version of the acceleration term in the geodesic equation and stress is a derivative of one vector field with respect to the other, half of a Lie bracket, defined in what are essentially local coordinates. For the landmark space case, we proved this formula in our previous paper \cite{MMM1}. We hope that the terms in this formula will be elucidated by further study and analysis of specific cases.

\section{\label{nmb:1} A Covariant Formula for Curvature}

\subsection{\label{nmb:1.1} Covariant derivative}
Let $(M,g)$ be a (finite dimensional) Riemannian manifold. 
There will be some formulae which are valid for infinite dimensional manifolds and we will introduce definitions for these below. For each $x\in M$ we view the metric also as a bijective mapping $g_x:T_xM\to T_x^*M$. Then $g\i$ is the metric on the cotangent bundle as well as the morphism $T^*M\to TM$. For a 1-form $\al\in \Om^1(M)= \Ga(T^*M)$ we consider the `sharp' vector field $\al^\sharp=g\i\al\in\X(M)$. If $\al = \al_i dx^i$, then $\al^\sharp = \al_i g^{ij} \partial_j$ 
is just the vector field obtained from $\al$ by `raising indices'. Similarly, for a vector field $X\in\X(M)$ we consider the `flat' 1-form $X^\flat = gX$. 
If $X=X^i \partial_i$, then $X^\flat = X^i g_{ij} dx^j$ 
is the 1-form obtained from $X$ by `lowering indices'. Note that
\begin{equation*}
\al(\be^{\sharp}) = g\i(\al,\be) = g(\al^\sharp,\be^\sharp) =\be(\al^{\sharp}).
\tag{1}\end{equation*}
Our aim is to express the sectional curvature of $g$ in terms of $\al,\be$ alone. It is important that the exterior derivative satisfies:
\begin{equation*}
d\al(\be^{\sharp},\ga^{\sharp}) 
= (\be^{\sharp})\al(\ga^{\sharp}) - (\ga^{\sharp})\al(\be^{\sharp}) -
\al([\be^{\sharp},\ga^{\sharp}])
\tag{2}\end{equation*}
We recall the definition of the Levi-Civita covariant derivative $\nabla$ and its basic properties:
\begin{align*}
2g(\nabla_XY,Z) &= X(g(Y,Z)) + Y(g(Z,X)) -Z(g(X,Y)) \tag{3}\\&\quad 
- g(X,[Y,Z])+ g(Y,[Z,X])+ g(Z,[X,Y])  \\
(\nabla_X\al)(Y) &= X\al(Y)-\al(\nabla_XY) \tag{4}\\
g((\nabla_X\al)^\sharp,Y)&= (\nabla_X\al)(Y)=X(\al(Y)) -\al(\nabla_XY) \\
& = Xg(\al^\sharp,Y) - g(\al^\sharp,\nabla_XY) = g(\nabla_X\al^\sharp,Y)
\quad \implies \\
\nabla_X(\al^\sharp) &= (\nabla_X\al)^\sharp \tag{5}\\
Xg\i(\al,\be)&=g\i(\nabla_X\al,\be)+g\i(\al,\nabla_X\be) \\
\nabla_{\al^{\sharp}}\be - \nabla_{\be^\sharp}\al &= g[\al^{\sharp},\be^{\sharp}] 
=[\al^{\sharp},\be^{\sharp}]^{\flat}
\end{align*}
From this follows
\begin{align*}
2(\nabla_{\al^{\sharp}}\be)(\ga^{\sharp}) &=
2g\i(\nabla_{\al^{\sharp}}\be,\ga) =
2g((\nabla_{\al^\sharp}\be)^\sharp,\ga^\sharp) =
2g(\nabla_{\al^\sharp}\be^\sharp,\ga^\sharp) =
\\&
= \al^{\sharp}g\i(\be,\ga) + \be^{\sharp}g\i(\ga,\al) -\ga^{\sharp}g\i(\al,\be)) 
\\&\quad
- g\i(\al,[\be^{\sharp},\ga^{\sharp}]^{\flat})
  +g\i(\be,[\ga^{\sharp},\al^{\sharp}]^{\flat})
  +g\i(\ga,[\al^{\sharp},\be^{\sharp}]^{\flat})
\\&
= \al^{\sharp}\be(\ga^{\sharp}) 
  +\be^{\sharp}\ga(\al^{\sharp}) 
  -\ga^{\sharp}\be(\al^{\sharp}) 
\\&\quad
- \al([\be^{\sharp},\ga^{\sharp}])
  +\be([\ga^{\sharp},\al^{\sharp}])
  +\ga([\al^{\sharp},\be^{\sharp}])
\\&
= \be^{\sharp}\ga(\al^{\sharp}) 
  -\al([\be^{\sharp},\ga^{\sharp}])
  +\ga([\al^{\sharp},\be^{\sharp}]) 
  - d\be(\ga^{\sharp},\al^{\sharp})
\\&
= \al^{\sharp}\ga(\be^{\sharp}) 
  -\be^{\sharp}\al(\ga^{\sharp})
  +\ga^{\sharp}\al(\be^{\sharp})
  +d\al(\be^{\sharp},\ga^{\sharp})
  -d\be(\ga^{\sharp},\al^{\sharp})
  -d\ga(\al^{\sharp},\be^{\sharp})
\tag{7}\end{align*}

\begin{theorem}\label{nmb:1.2} {\bf (Mario's Formula)}
Assume that all 1-forms $\al,\be,\ga,\de\in\Om_g^1(M)$ are closed. Then curvature is given by:
$$\boxed{
\begin{aligned}
&g\big(R(\al^{\sharp},\be^{\sharp})\ga^{\sharp},\de^{\sharp}\big) = R_1+R_2+R_3 \\
&\qquad R_1 = \tfrac14 \left(
-\al^{\sharp}\ga^{\sharp}\de(\be^{\sharp})
+\al^{\sharp}\de^{\sharp}\be(\ga^{\sharp})
+\be^{\sharp}\ga^{\sharp}\de(\al^{\sharp})
-\be^{\sharp}\de^{\sharp}\al(\ga^{\sharp})\right. \\
&\qquad\qquad  \left.
-\ga^{\sharp}\al^{\sharp}\de(\be^{\sharp})
+\ga^{\sharp}\be^{\sharp}\de(\al^{\sharp})
+\de^{\sharp}\al^{\sharp}\be(\ga^{\sharp})
-\de^{\sharp}\be^{\sharp}\al(\ga^{\sharp}) \right) \\
&\qquad R_2 = \tfrac14 \left(
-g\i\big(d(\ga(\be^{\sharp})),d(\de(\al^{\sharp}))\big)
+g\i\big(d(\ga(\al^{\sharp})),d(\de(\be^{\sharp}))\big) \right) \\
&\qquad R_3 = \tfrac14 \left(
g\big([\de^{\sharp},\al^{\sharp}],[\be^\sharp,\ga^{\sharp}]\big)
-g\big([\de^{\sharp},\be^{\sharp}],[\al^\sharp,\ga^{\sharp}]\big)
+2g\big([\al^{\sharp},\be^{\sharp}],[\ga^{\sharp},\de^{\sharp}]\big) \right)
\end{aligned}
}$$
For the numerator of sectional curvature we get
$$\boxed{
\begin{aligned}
&g\big(R(\al^{\sharp},\be^{\sharp})\be^{\sharp},\al^{\sharp}\big) = R_1+R_2+R_3 \\ 
&\qquad R_1 = \tfrac12 \left(
\al^{\sharp}\al^{\sharp}(\|\be\|^2)
-(\al^{\sharp}\be^{\sharp}+\be^{\sharp}\al^{\sharp})g\i(\al,\be)
+\be^{\sharp}\be^{\sharp}(\|\al\|^2)
 \right) 
\\&\qquad\quad
= \tfrac12\left(\al^\sharp\be([\al^\sharp,\be^\sharp])-\be^\sharp\al([\al^\sharp,\be^\sharp])\right)
\\&\qquad 
R_2 = \tfrac14 \left(
\|d(g\i(\al,\be))\|^2
-g\i\big(d(\|\al\|^2),d(\|\be\|^2) \right) \\
&\qquad R_3 = -\tfrac34
\big\|[\al^{\sharp},\be^{\sharp}]\big\|_g^2
\end{aligned}
}$$
\end{theorem}

Recall that sectional curvature is then
$$
k(\al^{\sharp},\be^{\sharp}) = 
\frac{g\big(R(\al^{\sharp},\be^{\sharp})\be^{\sharp},\al^{\sharp}\big)}{\|\al\|^2\|\be\|^2-g\i(\al,\be)^2}
$$

\begin{demo}{Proof}
We shall need that for a function $f$ we have:
\begin{align*}
(\nabla_{\be^{\sharp}}\ga)^{\sharp}f
&= df((\nabla_{\be^{\sharp}}\ga)^{\sharp}) 
= g\i(df,\nabla_{\be^{\sharp}}\ga)
= \be^{\sharp}g\i(df,\ga) - g\i(\nabla_{\be^{\sharp}}df,\ga)
\\&
= \be^{\sharp}\ga^{\sharp}f - \nabla_{\be^{\sharp}}df(\ga^{\sharp})
= \be^{\sharp}\ga^{\sharp}f 
  -\tfrac12\be^{\sharp}df(\ga^{\sharp}) 
  +\tfrac12 df^{\sharp}\ga(\be^{\sharp})
  -\tfrac12\ga^{\sharp}\be(df^{\sharp})
\\&
= \tfrac12 df^{\sharp}\ga(\be^{\sharp})
  +\tfrac12[\be^{\sharp},\ga^{\sharp}]f
= \tfrac12 d(\ga(\be^{\sharp}))(df^{\sharp})
  +\tfrac12[\be^{\sharp},\ga^{\sharp}]f
\tag{8}\end{align*}
For the three summands in the curvature formula, by multiple uses 
of formulas (2) and (7) and the closedness of $\al,\be,\ga,\de$, a straightforward calculation gives us:
\begin{align*}
&4(\nabla_{\al^{\sharp}}\nabla_{\be^{\sharp}}\ga)(\de^{\sharp}) =
\\&\qquad=
 2\al^{\sharp}(\nabla_{\be^{\sharp}}\ga)(\de^{\sharp})
-2(\nabla_{\be^{\sharp}}\ga)^{\sharp}\de(\al^{\sharp})
+2\de^{\sharp}\al((\nabla_{\be^{\sharp}}\ga)^{\sharp})
-2d(\nabla_{\be^{\sharp}}\ga)(\de^{\sharp},\al^{\sharp})
\\&\qquad= 
2\al^{\sharp}(\nabla_{\be^{\sharp}}\ga)(\de^{\sharp})
-d(\ga(\be^{\sharp}))(d(\de(\al^{\sharp})^{\sharp})
-[\be^{\sharp},\ga^{\sharp}]\de(\al^{\sharp})
\\&\qquad\quad
+2\al^{\sharp}(\nabla_{\be^{\sharp}}\ga)(\de^{\sharp})
+2(\nabla_{\be^{\sharp}}\ga)([\de^{\sharp},\al^{\sharp}])
= \cdots =
\\&\qquad=
-g\i\big(d(\ga(\be^{\sharp})),d(\de(\al^{\sharp}))\big)
+g\big([\de^{\sharp},\al^{\sharp}],[\be^\sharp,\ga^{\sharp}]\big)
\\&\qquad\quad
+2\al^{\sharp}\be^{\sharp}\ga(\de^{\sharp})
-2\al^{\sharp}\ga^{\sharp}\de(\be^{\sharp})
+\al^{\sharp}\de^{\sharp}\be(\ga^{\sharp})
-[\be^{\sharp},\ga^{\sharp}]\de(\al^{\sharp})
+\de^{\sharp}\al^{\sharp}\be(\ga^{\sharp})
\\ &\text{and similarly}
\\ &
-4(\nabla_{\be^{\sharp}}\nabla_{\al^{\sharp}}\ga)(\de^{\sharp})
=
+g\i\big(d(\ga(\al^{\sharp})),d(\de(\be^{\sharp}))\big)
-g\big([\de^{\sharp},\be^{\sharp}],[\al^\sharp,\ga^{\sharp}]\big)
\\&\qquad\quad
-2\be^{\sharp}\al^{\sharp}\ga(\de^{\sharp})
+2\be^{\sharp}\ga^{\sharp}\de(\al^{\sharp})
-\be^{\sharp}\de^{\sharp}\al(\ga^{\sharp})
+[\al^{\sharp},\ga^{\sharp}]\de(\be^{\sharp})
-\de^{\sharp}\be^{\sharp}\al(\ga^{\sharp})
\\&
-2(\nabla_{[\al^{\sharp},\be^{\sharp}]}\ga)(\de^{\sharp}) =
\\&\qquad=
-[\al^{\sharp},\be^{\sharp}]\ga(\de^{\sharp})
+\ga^{\sharp}\de([\al^{\sharp},\be^{\sharp}])
-\de^{\sharp}\ga([\al^{\sharp},\be^{\sharp}])
-d[\al^{\sharp},\be^{\sharp}]^{\flat}(\ga^{\sharp},\de^{\sharp})
\\&\qquad=
-[\al^{\sharp},\be^{\sharp}]\ga(\de^{\sharp})
+g\big([\al^{\sharp},\be^{\sharp}],[\ga^{\sharp},\de^{\sharp}]\big)
\end{align*}
Now we can compute the curvature (remember that $d\al=d\be=\dots =0$):
\begin{align*}
&4g\big(R(\al^{\sharp},\be^{\sharp})\ga^{\sharp},\de^{\sharp}\big)
=4\de\big(R(\al^{\sharp},\be^{\sharp})\ga^{\sharp}\big)
=4\de\big(\nabla_{\al^{\sharp}}\nabla_{\be^{\sharp}}\ga^{\sharp}
-\nabla_{\be^{\sharp}}\nabla_{\al^{\sharp}}\ga^{\sharp}
-\nabla_{[\al^{\sharp},\be^{\sharp}]}\ga^{\sharp}\big)
\\&\qquad=
4\big(\nabla_{\al^{\sharp}}\nabla_{\be^{\sharp}}\ga
-\nabla_{\be^{\sharp}}\nabla_{\al^{\sharp}}\ga
-\nabla_{[\al^{\sharp},\be^{\sharp}]}\ga\big)(\de^\sharp)
\\&\qquad=
-g\i\big(d(\ga(\be^{\sharp})),d(\de(\al^{\sharp}))\big)
+g\i\big(d(\ga(\al^{\sharp})),d(\de(\be^{\sharp}))\big)
\\&\qquad\quad
+g\big([\de^{\sharp},\al^{\sharp}],[\be^\sharp,\ga^{\sharp}]\big)
-g\big([\de^{\sharp},\be^{\sharp}],[\al^\sharp,\ga^{\sharp}]\big)
+2g\big([\al^{\sharp},\be^{\sharp}],[\ga^{\sharp},\de^{\sharp}]\big)
\\&\qquad\quad
-\al^{\sharp}\ga^{\sharp}\de(\be^{\sharp})
+\al^{\sharp}\de^{\sharp}\be(\ga^{\sharp})
+\be^{\sharp}\ga^{\sharp}\de(\al^{\sharp})
-\be^{\sharp}\de^{\sharp}\al(\ga^{\sharp})
\\&\qquad\quad
-\ga^{\sharp}\al^{\sharp}\de(\be^{\sharp})
+\ga^{\sharp}\be^{\sharp}\de(\al^{\sharp})
+\de^{\sharp}\al^{\sharp}\be(\ga^{\sharp})
-\de^{\sharp}\be^{\sharp}\al(\ga^{\sharp})
\end{align*}
For the sectional curvature expression this simplifies 
(as always, for closed 1-forms) to the expression in the theorem. 
The two versions of $R_1$ correspond to each other, using $d\al=0$ and $d\be=0$. 
\qed \end{demo}

\subsection{\label{nmb:1.3} Mario's formula in coordinates}
The   formula   for   sectional  curvature  becomes  especially
transparent  if we expand it in coordinates. Assume that $\al =
\al_i  dx^i,  \be  = \be_i dx^i$ where the coefficients $\al_i,
\be_i$  are  {\it constants}, hence $\al, \be$ are closed. Then
$\al^\sharp  = g^{ij}\al_i \partial_j, \be^\sharp = g^{ij}\be_i
\partial_j$.  Substituting these in the terms of the right hand
side of Mario's formula for sectional curvature, we get:
\begin{align*}
&\text{2nd deriv. terms} = 2R_1 
=  2\al^{\sharp}\al^{\sharp}(\|\be\|^2)
+2\be^{\sharp}\be^{\sharp}(\|\al\|^2)
-2(\al^{\sharp}\be^{\sharp}+\be^{\sharp}\al^{\sharp})g\i(\al,\be)
\\&\qquad = 2\al_i g^{is}(\al_j g^{jt}( \be_k \be_l g^{kl})_{,t} )_{,s}
+ 2\be_i g^{is}(\be_j g^{jt}( \al_k \al_l g^{kl})_{,t} )_{,s}
\\& \qquad - 2\al_i g^{is}(\be_j g^{jt}( \be_k \al_l g^{kl})_{,t} )_{,s}
 - 2\be_i g^{is}(\al_j g^{jt}( \al_k \be_l g^{kl})_{,t} )_{,s}
\\&\qquad = 2(\al_i \be_k - \al_k \be_i)\cdot 
(\al_j \be_l - \al_l \be_j)\cdot g^{is}(g^{jt} g^{kl}_{,t})_{,s}
\\& \text{1st deriv. terms} = 4 R_2 = \|d(g\i(\al,\be))\|^2
-g\i\big(d(\|\al\|^2),d(\|\be\|^2)\big)
\\& \qquad = (\al_i \be_j g^{ij})_{,s} g^{st} (\al_l \be_k g^{kl})_{,t}
- (\al_i \al_j g^{ij})_{,s} g^{st} ( \be_k \be_l g^{kl})_{,t}
\\& \qquad = -\tfrac12 (\al_i \be_k - \al_k \be_i)\cdot 
(\al_j \be_l - \al_l \be_j)\cdot g^{ij}_{,s}g^{st} g^{kl}_{,t}
\\& \text{Lie bracket} = [\al^{\sharp},\be^{\sharp}] 
= \left(\al_i g^{is}(\be_k g^{kt})_{,s}-\be_i g^{is}(\al_k g^{kt})_{,s}\right) \partial_t
\\&\qquad = (\al_i \be_k - \al_k \be_i) g^{is} g^{kt}_{,s}\partial_t
\\& \text{Lie bracket term} = 4R_3 
= -3g\big([\al^{\sharp},\be^{\sharp}],[\al^{\sharp},\be^{\sharp}]\big)
\\& \qquad =-3 (\al_i \be_k - \al_k \be_i)\cdot 
(\al_j \be_l - \al_l \be_j)\cdot g^{is}g^{kp}_{,s} g_{pq} g^{jt} g^{lq}_{,t}
\end{align*}
hence we have the coordinate version for the three terms in sectional curvature:
$$\boxed{
\begin{aligned}
&g\big(R(\al^{\sharp},\be^{\sharp})\be^{\sharp},\al^{\sharp}\big) 
= (\al_i \be_k - \al_k \be_i)\cdot (\al_j \be_l - \al_l \be_j) \left( R_1^{ijkl} + R_2^{ijkl} + R_3^{ijkl} \right) 
\\ 
&\qquad  R_1^{ijkl} = \tfrac12\cdot g^{is}(g^{jt} g^{kl}_{,t})_{,s} 
\\
&\qquad R_2^{ijkl} = -\tfrac18\cdot g^{ij}_{,s}g^{st} g^{kl}_{,t} 
\\
&\qquad R_3^{ijkl} = -\tfrac34\cdot g^{is}g^{kp}_{,s} g_{pq} g^{jt} g^{lq}_{,t}
\end{aligned}
}$$
Note that the usual contravariant metric tensor $g_{ij}$ occurs
in  only  one  place,  everything  else  being derived from the
covariant metric tensor $g^{ij}$. 
Note that the fist term $R_1$ can be split into a pure second derivative term 
$R_{11}=g^{is}g^{jt}g^{kl}_{,st}$ plus a first derivative term $R_{12}=g^{is}g^{jt}_{,s} g^{kl}_{,t}$.

There is also a version of Mario's formula which is, in a sense, 
intermediate between the coordinate free version and the coordinate version. 
The main thing that coordinates allow you to do is to take derivatives using the associated flat connection. 
In the case of this formula, this introduces auxiliary vector fields $X_\al$ and 
$X_\be$ playing the role of `locally constant' extensions of the value of $\al^\sharp$ and 
$\be^\sharp$ at the point $x \in M$ where the curvature is being calculated and for which 
the 1-forms $\al, \be$ appear locally constant too. More precisely, assume we are given $X_\al$ and $X_\be$ such that:
\begin{enumerate}
\item $X_\al(x) = \al^\sharp(x),\quad X_\be(x) = \be^\sharp(x)$,
\item Then $\al^\sharp-X_\al$ is zero at $x$ hence has a well defined derivative 
       $D_x(\al^\sharp-X_\al)$ lying in Hom$(T_xM,T_xM)$. For a vector field $Y$ we have 
       $D_x(\al^\sharp-X_\al).Y_x = [Y,\al^\sharp-X_\al](x) = \L_Y(\al^\sharp-X_\al)|_x$.
       The same holds for $\be$.
\item $\L_{X_\al}(\al)=\L_{X_\al}(\be)=\L_{X_\be}(\al)=\L_{X_\be}(\be)=0$,
\item $[X_\al, X_\be] = 0$.
\end{enumerate}
Locally constant 1-forms and vector fields satisfy these properties. Using these forms and vector fields, we then define:
\begin{align*}
\mathcal F(\al,\be) :&= \tfrac12 d(g\i(\al,\be)), \qquad \text{a 1-form on $M$ called the {\it force},}\\
\mathcal D(\al,\be)(x) :&= D_x(\be^\sharp - X_\be).\al^\sharp(x)
\\&
= d(\be^\sharp - X_\be).\al^\sharp(x), \quad\text{a tangent vector at $x$ called the {\it stress}.}
\end{align*}  
Then in the notation above:
$$\boxed{ 
\begin{aligned}
&g\big(R(\al^{\sharp},\be^{\sharp})\be^{\sharp},\al^{\sharp}\big)(x) = R_{11} +R_{12} + R_2 + R_3 \\ 
&\quad R_{11} = \tfrac12 \left(
\L_{X_\al}^2(g\i)(\be,\be)-2\L_{X_\al}\L_{X_\be}(g\i)(\al,\be)
+\L_{X_\be}^2(g\i)(\al,\al)
 \right)(x) \\
& \quad R_{12} = \langle \mathcal F(\al,\al), \mathcal D(\be,\be) \rangle + \langle \mathcal F(\be,\be),\mathcal D(\al,\al)\rangle  
- \langle \mathcal F(\al,\be), \mathcal D(\al,\be)+\mathcal D(\be,\al) \rangle \\
&\quad R_2 = \left(
\|\mathcal F(\al,\be)\|^2_{g\i}
-\big\langle \mathcal F(\al,\al)),\mathcal F(\be,\be)\big\rangle_{g\i} \right)(x) \\
&\quad R_3 = -\tfrac34 \| \mathcal D(\al,\be)-\mathcal D(\be,\al) \|^2_{g_x} \label{eq:curv3}
\end{aligned}
}$$

The reformulation of $R_1$ follows from the calculation:
\begin{align*}
\al^\sharp \al^\sharp(\|\be\|^2)(x) &= X_\al \al^\sharp(\|\be\|^2)(x) \\ 
&=X_\al X_\al(\|\be\|^2)(x) + X_\al (\al^\sharp - X_\al)(\|\be\|^2)(x) \\
&= \L_{X_\al} \L_{X_\al}(g\i(\be, \be))(x) + \langle D_x(\al^\sharp - X_\al).X_\al(x),d\|\be\|^2)(x)\rangle \\
&= \L_{X_\al}^2(g\i)(\be, \be)(x) + \langle D_x(\al^\sharp - X_\al).X_\al(x),d\|\be\|^2)(x)\rangle
\end{align*}
and the similar result for the other terms. The reformulation of $R_3$ comes from the calculation:
\begin{align*}
[\al^\sharp, \be^\sharp](x) &= (X_\al \circ \be^\sharp)(x) - (X_\be \circ \al^\sharp)(x) \\
&= (X_\al \circ (\be^\sharp - X_\be))(x) - (X_\be \circ (\al^\sharp - X_\al)(x)\\
&= D_x((\be^\sharp - X_\be).X_\al(x) - D_x(\al^\sharp - X_\al).X_\be(x)
\end{align*}

\subsection{\nmb.{1.4} Infinite dimensional manifolds}
The main focus of this paper are the infinite dimensional manifolds of diffeomorphisms of a finite dimensional $N$, of the embeddings of one finite dimensional $M$ into another $N$ and of the set of submanifolds $F$ of a manifold $N$. These are infinite dimensional and can be realized in multiple ways depending on the degree of smoothness imposed on the diffeomorphism/embedding/submanifold. The first two have realizations as Hilbert manifolds but the last does not. Moreover, the group law on the Hilbert manifold version of the group of diffeomorphisms is not differentiable. If one desires to carry over finite dimensonal techniques to the infinite dimensional setting, it works much more smoothly to use the Frechet space of $C^\infty$ functions decreasing rapidly at infinity as the base vector space for charts of these spaces. But then its dual is not Frechet, so one needs a bigger category for charts on bundles. The best setting has been developed by one of the authors and his collaborators \cite{KM97} and uses `$c^\infty$-open' subsets in arbitrary `convenient' locally convex topological vector spaces for charts. This theory and some of the reasons why it works are summarized in the appendix. For our purposes, complete locally convex topological vector spaces (which are always convenient) suffice and, on them `$c^\infty$-open' just means open.

To extend Mario's formula to infinite-dimensional manifolds then, let  $(M,g)$  be  a so-called `weak Riemannian manifold' \cite{KM97}: a convenient manifold $M$ and smooth map:
$$g: TM \times_M TM \longrightarrow \R$$
which is a positive definite symmetric bilinear form $g_x$ on each tangent space $T_xM, x \in M$. For a convenient  manifold  we have to choose what we mean by 1-forms carefully. For each $x\in M$ the metric defines a  mapping  $g_x:T_xM\to  T_x^*M$ (which we denote by the same symbol $g_x$). In the case of a Riemannian Hilbert manifold, this is bijective and has an inverse but otherwise is only injective, hence the term `weak metric'. The image $g(TM)\subset  T^*M$  is  called  the {\it g-smooth cotangent bundle}.  Then  $g\i$ is the metric on the $g$-smooth cotangent bundle as well as the morphism $g(TM)\to TM$. Now  define $\Om_g^1(M):=  \Ga(g(TM))$ and 
$\al^\sharp=g\i\al\in\X(M), X^\flat  =  gX$  are  as  above. 
The exterior derivative is now defined by:
\begin{equation*}
d\al(\be^{\sharp},\ga^{\sharp}) 
= (\be^{\sharp})\al(\ga^{\sharp}) - (\ga^{\sharp})\al(\be^{\sharp}) -
\al([\be^{\sharp},\ga^{\sharp}])
\end{equation*}
We have $d:\Om_g^1(M)\to \Om^2(M)=\Ga(L^2_{\text{skew}}(TM;\mathbb R))$ 
since the embedding $g(TM)\subset T^*M$ is a smooth fiber linear mapping. 
Note  that  on  an  infinite dimensional manifold $M$ there are
many  choices  of  differential  forms  but only one of them is
suitable  for  analysis  on  manifolds.  These are discussed in
\cite[Section  33]{KM97}.  Here  we  consider  subspaces  of these
differential forms.

Further requirements need to be imposed on $(M,g)$ for our theory to work. Since it is an infinite 
dimensional weak   Riemannian   manifold  the Levi-Civita  covariant  derivative  might  not  exist  in $TM$. 
The Levi-Civita covariant derivative exists if and only if  
the metric itself admits gradients with respect to itself in the following senses. 
The easiest way to express this is locally in a chart $U \subset M$. Let $V_U$ be the vector space 
of constant vector fields on $U$. 
Then we assume that there are smooth maps grad$_1 g$ and grad$_2 
g$ from $U\times V_U$ to $V_U$, quadratic in $V_U$  such that
\begin{equation*} \boxed{
\begin{aligned}
D_{x,Z}g_x(X,X) &= g_x(Z,\on{grad}_1g(x)(X,X))
\\ 
D_{x,X}g_x(X,Z) &= g_x(\on{grad}_2g(x)(X,X),Z)
\end{aligned}\;\text{  for all }Z.
}\end{equation*}
(If we express this globally we also get derivatives of the vector fields
$X$ and $Z$.) 
This allows to use \thetag{\nmb!{1.1}.3} to get the covariant derivative. 
Then the rest of the derivation of Mario's formula goes through and  the  final  formula for curvature holds in both the finite and infinite dimensional cases.
There are situations where the covariant derivative exists but not both gradients; see 
\cite{BauerMichor2012}, and the corresponding extension of \cite[appendix]{YMSM2008} to the real 
line.

Some  constructions to be done shortly encounter a second problem: they lead to vector fields whose values do not lie in $T_x M$, but in the Hilbert space completion $\overline{T_x M}$ with respect to the inner product $g_x$. To manipulate these as in the finite dimensional case, we need to know that $\bigcup_{x\in M}\overline{T_x M}$ forms a smooth vector bundle over $M$. More precisely, choose an atlas 
$(U_\al,u_\al : U_\al\to E)$ of $M$, where the $U_\al\subset M$ form an open cover of $M$, where each $u_\al:U_\al\to E$ is a homeomorphism of $U_\al$ onto the open subset $u_\al(U_\al)$ of the convenient vector space $E$ which models $M$, and where $u_{\al\be}= u_\al\o u_\be\i: u_\be(U_\al\cap U_\be)\to u_\al(U_\al\cap U_\be))$ is a smooth diffeomorphism. The mappings 
$x\mapsto \ph_{\al\be}(x) = d u_{\al\be}(u_\be\i(x))\in L(E,E)$ then form the cocycle of transition functions $\ph_{\al\be}:U_\al\cap U_\be \to GL(E)$ wich define the tangent bundle $TM$. We then assume that the local expression of each Riemannian metric $g_x$ on $E$ are equivalent weak inner products hence define Hilbert space completions which are quasi-isometric via extensions of the embeddings of $E$ (in each chart). Let us call one such Hilbert space $\mathcal H$. 
We then require that all transition functions $\ph_{\al\be}(x):E\to E$ extend to bounded linear isomorphisms $\mathcal H\to \mathcal H$ and that each $\ph_{\al\be}:U_\al\cap U_\be\to L(\mathcal H,\mathcal H)$ is again smooth.

These two properties will be sufficient
for  all  the  constructions  we  need  so  we make them into a
definition:

\subsection*{Definition}  A convenient weak Riemannian manifold
$(M,g)$ will be called a {\it robust} Riemannian manifold if:
\begin{enumerate}
  \item The Levi-Civita covariant derivative exists. Equivalently, 
	the metric $g_x$ admits gradients in the above two senses.
  \item The completions $\overline{T_x M}$ form a vector bundle as described above.
\end{enumerate}

Note that a Hilbert manifold is automatically robust. We can make the relationship between robust 
manifolds and Hilbert manifolds more explicit if we introduce another definition, that of a 
{\it pre Hilbert manifold} similar to the notion of a pre-Hilbert topological vector space:

\subsection*{Definition} A robust Riemannian manifold $(M,g)$ is called {\it pre-Hilbert} if there 
exists an atlas $(U_\al,u_\al : U_\al\to E)$ for which:
\begin{enumerate}
  \item Each $u_\al(U_\al)$ is contained in the Hilbert norm interior of its closure in 
         $\mathcal H$, which we denote $u_\al(U_\al)^{\mathcal H}$. 
  \item All chart change maps $u_{\al\be}$ extend to smooth mappings between the open subsets 
         $u_\al(U_\al)^{\mathcal H}$, hence define a completion $M \subset M^{\mathcal H}$ which is 
         a Hilbert manifold.
\end{enumerate}

Note that in this definition the atlas must be properly chosen: for example its open sets $U_\al$ 
must be open in the weak topology defined by path lengths. More precisely, for any weak Riemannian 
manifold $M$, the inner products $g_x$ assign a length to every smooth path in $M$ and we get a 
distance function $d(x,y)$ as the infimum of lengths of paths joining $x$ and $y$ (which might 
however be zero for some $x \ne y$). The topology defined by path lengths is usually much weaker 
than the strong topology given by the definition of $M$.

These distinctions are well illustrated by the spaces we will discuss below. Firstly, manifolds of 
smooth mappings like $\on{Emb}(M,N)$ with their canonically induced Sobolev metrics of order 
$s>\dim{M}/2$ do admit completions  $\on{Emb}^s(M,N)$ to Hilbert manifolds hence are pre-Hilbert; 
see \cite[42.1]{KM97} for the explicit chart changes.   
But their quotient manifolds $B(M,N)=\on{Emb}(M,N)/\on{Diff}(M)$ are only robust in general because 
the second condition fails. The extensions of the chart change maps are homeomorphisms but not 
differentiable: This is due to the fact that the Sobolev completions $\on{Diff}^s(M)$ 
of $\on{Diff}(M)$ of order $s>\dim{M}/2$ are smooth manifolds themselves, but only topological 
groups: right translations are still smooth, left translations and inversions however, are only 
continuous (and not even Lipschitz). So the action of $\on{Diff}^s(M)$ on $\on{Emb}^s(M,N)$, after 
Sobolev completion, has aspects which are only continuous and thus $B^s(M,N)= 
\on{Emb}^s(M,N)/\on{Diff}^s(M)$ is only a topological manifold in general. This phenomenon also 
appears in the chart changes of the canonical atlas of $B(M,N)$; see \cite[44.1]{KM97} 
for an explicit formula of the chart change and the role of inversion in $\on{Diff}^s(M)$ in it.        

\subsection{\label{nmb:1.5} Covariant curvature and O'Neill's formula, finite dimensional}
Let  $p:(E,g_E)\to  (B,g_B)$ be a Riemannian submersion between finite dimensional manifolds, 
i.e., for each $b\in B$ and $x\in E_b:=p\i(b)$ the 
$g_E$-orthogonal splitting $T_xE =  T_x(E_{p(x)})\oplus T_x(E_{p(x)})^\bot 
=:T_x(E_{p(x)})\oplus \on{Hor}_x(p)$ has the property that 
$T_xp:(\on{Hor}_x(p), g_E)\to (T_bB,g_B)$ is an isometry. Each vector field $X\in\X(E)$ 
is decomposed as $X=X^{\text{hor}}+X^{\text{ver}}$ into horizontal and vertical parts. 
Each vector field $\xi\in\X(B)$ can be uniquely lifted to a smooth horizontal field 
$\xi^{\text{hor}}\in \Ga(\on{Hor}(p))\subset \X(E)$. 
O'Neill's formula says that for any two horizontal vector fields $X,Y$ on $E$ 
and any $x \in E$, the sectional curvatures of $E$ and $B$ are related by:
$$ g_{p(x)}(R^B(p_*(X_x),p_*(Y_x))p_*(Y_x),p_*(X_x)) = g_x(R^E(X_x,Y_x)Y_x,X_x) + \tfrac34 \|[X,Y]^{ver}\|^2_x.$$
Comparing Mario's formula on $E$ and $B$ gives an immediate proof of this fact. Start with:

\begin{lemma*}
If $\al\in\Om^1(B)$ is a 1-form on $B$, then the vector field $(p^*\al)^\sharp$ 
is horizontal and we have $Tp\o (p^*\al)^\sharp = \al^\sharp\o p$.
Therefore $(p^*\al)^\sharp$ equals the horizontal lift $(\al^\sharp)^{\text{hor}}$. 
For each $x\in E$ the mapping $(T_xp)^*:(T_{p(x)}^*B, g_B\i)\to (T_x^*E,g_E\i)$ is an isometry.
\end{lemma*}

\begin{demo}{Proof}
All this holds because for $X_x\in T_xE$ we have:
\begin{align*}
g_E((p^*\al)^\sharp_x,X_x)&= (p^*\al)_x(X_x)= \al_{p(x)}(T_xp.X_x)
=\al_{p(x)}(T_xp.X_x^{\text{hor}})
\\&
=g_E((p^*\al)^\sharp_x,X_x^{\text{hor}}) \\
g_B(T_xp(p^*\al)^\sharp_x,T_xp.X_x) &= g_E((p^*\al)^\sharp_x,X_x^{\text{hor}})
= \al_{p(x)}(T_xp.X_x) 
\\&
= g_B(\al_{p(x)}^\sharp,T_xp.X_x). \qed
\end{align*}
\end{demo}

More generally we have:
$$ 
g_E\i(p^*\al,p^*\be) = g_E((p^*\al)^\sharp,(p^*\be)^\sharp) 
= g_B(\al^\sharp, \be^\sharp)\o p =p^* g_B\i(\al,\be).
$$
Consequently, we get for 1-forms $\al$, $\be$ on $B$:
\begin{align*}
d\|p^*\al\|_{g_E\i}^2 &= d p^*\|\al\|_{g_B\i}^2 = p^* d\|\al\|_{g_B\i}^2
\\
(p^*\be)^\sharp \|p^*\al\|_{g_E\i}^2 &= (p^* d\|\al\|_{g_B\i}^2)((\al^\sharp)^{\text{hor}})
  = p^*(\be^\sharp\|\al\|_{g_B\i}^2)
\end{align*}
In the following computation we use
\begin{align*}
\|[(p^*\al)^{\sharp},(p^*\be)^{\sharp}]^{\text{hor}}\|_{g_E}^2
&
=p^*\|[\al^{\sharp},\be^{\sharp}]\|_{g_B}^2
\end{align*}
We take Mario's formula \thetag{\ref{nmb:1.2}} and apply it to the closed 1-forms 
$p^*\al$, $p^*\be$ on $E$ where $\al$, $\be$ are closed 1-forms on $B$. 
Using the results above we get:
\begin{align*}
&4g_E\big(R((p^*\al)^{\sharp},(p^*\be)^{\sharp})(p^*\be)^{\sharp},(p^*\al)^{\sharp}\big) =
\\&\qquad=
\|d(g_E\i(p^*\al,p^*\be))\|_{g_E\i}^2
-g_E\i\big(d(\|p^*\al\|_{g_E\i}^2),d(\|p^*\be\|_{g_E\i}^2)\big)
\\&\qquad\quad
-3\|[(p^*\al)^{\sharp},(p^*\be)^{\sharp}]^{\text{hor}}\|_{g_E}^2
-3\|[(p^*\al)^{\sharp},(p^*\be)^{\sharp}]^{\text{ver}}\|_{g_E}^2
\\&\qquad\quad
+2(p^*\al)^{\sharp}(p^*\al)^{\sharp}(\|p^*\be\|_{g_E\i}^2)
+2(p^*\be)^{\sharp}(p^*\be)^{\sharp}(\|p^*\al\|_{g_E\i}^2)
\\&\qquad\quad
-2((p^*\al)^{\sharp}(p^*\be)^{\sharp}+(p^*\be)^{\sharp}(p^*\al)^{\sharp})g_E\i(p^*\al,p^*\be)
\\&\qquad=
p^*\|d(g_B\i(\al,\be))\|_{g_B\i}^2
-p^* g_B\i\big(d(\|\al\|_{g_B\i}^2),d(\|\be\|_{g_B\i}^2)\big)
-3 p^*\|[\al^{\sharp},\be^{\sharp}]\|_{g_B}^2
\\&\qquad\quad
-3\|[(p^*\al)^{\sharp},(p^*\be)^{\sharp}]^{\text{ver}}\|_{g_E}^2
\\&\qquad\quad
+2p^*\big(\al^{\sharp}\al^{\sharp}(\|\be\|_{g_B\i}^2)\big)
+2p^*\big(\be^{\sharp}\be^{\sharp}(\|\al\|_{g_B\i}^2)\big)
-2p^*\big((\al^{\sharp}\be^{\sharp}+\be^{\sharp}\al^{\sharp})g_B\i(\al,\be)\big)
\\&\qquad
=4p^* g_B\big(R^B(\al^{\sharp},\be^{\sharp})\be^{\sharp},\al^{\sharp}\big) 
-3\|[(p^*\al)^{\sharp},(p^*\be)^{\sharp}]^{\text{ver}}\|_{g_E}^2
\end{align*}
which is a short proof of O'Neill's formula.

\subsection{\label{nmb:1.6} Covariant curvature and O'Neill's formula}
Let $p:(E,g_E)\to (B,g_B)$ be a Riemann submersion between infinite dimensional 
robust Riemann manifolds; i.e., for each $b\in B$ and $x\in E_b:=p\i(b)$ 
the tangent mapping $T_xp:(T_xE, g_E)\to (T_bB,g_B)$ 
is a surjective metric quotient map so that 
\begin{equation*}
\|\xi_b\|_{g_B} := \inf\bigl\{X_x\in T_xE: T_xp.X_x=\xi_b\bigr\}.
\tag{1}\end{equation*}
The infinimum need not be attained in $T_xE$ but will be in the completion 
$\overline{T_xE}$. 
The orthogonal subspace $\{Y_x: g_E(Y_x,T_x(E_b))=0\}$ 
has therefore to be taken in $\overline{T_xE}$.

If $\al_b=g_B(\al_b^\sharp,\quad)\in g_B(T_bB)\subset T_b^*B$ 
is an element in the $g_B$-smooth dual, then 
$p^*\al_b:=(T_xp)^*(\al_b)= g_B(\al_b^\sharp,T_xp\quad):T_xE\to \mathbb R$ 
is in $T_x^*M$ but in general it is not an element in the smooth dual $g_E(T_xE)$. 
It is, however, an element of the Hilbert space completion $\overline{g_E(T_xE)}$ 
of the $g_E$-smooth dual $g_E(T_xE)$ with respect to the norm $\|\quad\|_{g_E\i}$, 
and the element $g_E\i(p^*\al_b)=: (p^*\al_b)^\sharp$ is in the 
$\|\quad\|_{g_E}$-completion $\overline{T_xE}$ of $T_xE$. 
We can call $g_E\i(p^*\al_b)=: (p^*\al_b)^\sharp$ the {\it horizontal lift} of 
$\al_b^\sharp = g_B\i(\al_b)\in T_bB$.

In the following we discuss the manifold $E$ and we write $g$ instead of $g_E$. The metric $g_x$ can be evaluated at elements in the completion 
$\overline{T_xE}$. Moreover, for any smooth sections $X,Y\in\Ga(\overline{TE})$ the mapping $g(X,Y):M\to \mathbb R$ is still smooth: This is a local question, so let $E$ be $c^\infty$-open in a convenient vector space $V_E$.
Since the evaluations on $X\otimes Y$ form a set of bounded linear functionals on the space 
$L^2_{\text{sym}}(\overline{V_M};\mathbb R)$ of bounded symmetric bilinear forms on 
$\overline{V_M}$ which recognizes bounded subsets, it follows that $g$ is smooth as a mapping 
$M\to L^2_{\text{sym}}(\overline{V_M};\mathbb R)$, by the smooth uniform boundedness theorem, see \cite{KM97}.

\begin{lemma*}
If $\al$ is a smooth 1-form on an open subset $U$ of $B$ with values in the $g_B$-smooth dual $g_B(TB)$, 
then $p^*\al$ is a smooth 1-form on $p\i(U)\subset E$ with values in the $\|\quad\|_{g_E\i}$-completion 
of the $g_E$-smooth dual $g_E(TE)$. Thus also $(p^*\al)^\sharp$ is smooth from $E$ into the 
$g_E$-completion of $TE$, and it has values in the $g_E$-orthogonal subbundle to the vertical 
bundle in the $g_E$-completion.
We may continuously extend $T_xp$ to the $\|\quad\|_{g_E\i}$-completion, 
and then we have $Tp\o (p^*\al)^\sharp= \al^\sharp\o p$. 
Moreover, the Lie bracket of two such forms, $[(p^*\al)^\sharp,(p^*\be)^\sharp]$, is defined. 
The exterior derivative $d(p^*\al)$ is defined and is applicable to vector fields with values 
in the completion like $(p^*\be)^\sharp$. 
\end{lemma*}

That the Lie bracket is defined, is also a non-trivial statement: 
We have to differentiate in directions which are not tangent to the manifold. 

\begin{demo}{Proof of the lemma}
This is a local question; so we may assume that $U=B$ and $p\i(U)=E$  are $c^\infty$-open subsets in convenient 
vector spaces $V_B$ and $V_E$, respectively, so that all tangent bundles are trivial. 
By definition, $\al^\sharp = g_B\i\o \al:B\to B\x V_B$ is smooth. We have to 
show that $(p^*\al)^\sharp = g_E\i\o p^*\al$ is a smooth mapping from $E$ into the 
$\|\quad\|_{g_E}$-completion of $V_E$. By the smooth uniform boundedness theorem (see \cite{KM97})
it suffices to check that the composition with each bounded linear functional in a set 
$\mathcal S\subset V_E'$ is smooth, where $\mathcal S\subseteq V_E'$ is a set of linear functionals 
on $V_E$ which recognizes bounded subset of $V_E$. For this property, functionals of the form
$g_E(v,\quad)$ for $v\in V_E$ suffice. But 
$$x\mapsto (g_E)_x(v,(p^*\al)^\sharp|_x ) = p^*\al|_x(v)  = \al|_x(T_xp.v)$$
is obviously smooth. 

We may continuously extend the metric quotient mapping $T_xp$ to the 
$\|\quad\|_{g_E}$-completion and get a mapping $T_xp:\overline{T_xE}\to\overline{T_{b}B}$ where 
$b=p(x)$. 
For a second form $\be\in \Ga(g_B(TB))$ we have then 
\begin{multline*}
g_B(\be^\sharp|_b,T_xp.(p^*\al)^\sharp|_x)= (\be_b(T_xp.(p^*\al)^\sharp|_x) = 
(p^*\be)|_x((p^*\al)^\sharp|_x) =
\\
= g_E\i((T_xp)^*\be,(T_xp)^*\al)  = g_B(\be_b,\al_b) = g_B(\be^\sharp|_b,(\al^\sharp\o p)(x))
\end{multline*}
which implies $Tp\o (p^*\al)^\sharp = \al^\sharp\o p$.

For the Lie bracket of two such forms, $[(p^*\al)^\sharp,(p^*\be)^\sharp]$, we can again assume that 
all bundles are trivial. Then
\begin{align*}
[(p^*\al)^\sharp,(p^*\be)^\sharp](x) &= d((p^*\be)^\sharp)(x)((p^*\al)^\sharp)
- d((p^*\al)^\sharp)(x)((p^*\be)^\sharp)
\\
d((p^*\be)^\sharp)(x)((p^*\al)^\sharp) &= d(g_E\i\o (Tp)^*\o \be\o p)(x)((p^*\al)^\sharp)
\\&
= d(g_E\i\o (Tp)^*\o \be)(b).T_xp.(p^*\al)^\sharp
\\&
= d(g_E\i\o (Tp)^*\o \be)(b).\al^\sharp(p(x)).
\end{align*}
So the Lie bracket is well defined.
\qed\end{demo}

By assumption, the metric $g=g_E$ admits gradients with respect to itself as in 
\thetag{\ref{nmb:1.4}}. In a local chart we have
\begin{align*}
D_{x,Z}g_x(X,X) &= g_x(Z,\on{grad}_1g(x)(X,X))
\\ 
D_{x,Z}g_x(Z,X) &= g_x(\on{grad}_2g(x)(Z,Z),X)
\tag{2}\end{align*}
for $X,Z\in V_E$. 
We can then take $X,\in\overline{V_E}$ in the upper left expression of \thetag{2} and thus 
also in the right hand side. Then the upper right term allows to take $Z\in \overline{V_EM}$ also. 
This carries over to the lower expression.


Thus the local expressions of the Christoffel symbols of the Levi-Civita covariant derivative 
extend to sections of the completed tensor bundle 
$\overline{TE}$, and therefore the Levi-Civita covariant derivative extends to smooth sections of 
$\overline{TE}$ which are differentiable in directions in $\overline{TE}$ like $(p^*\al)^\sharp$.
Thus expressions like $\nabla^E_{(p^*\al)^\sharp}(p^*\be)^\sharp$ make sense and are again of the same 
type so that one can iterate. Thus the curvature expression 
$g_E\big(R((p^*\al)^{\sharp},(p^*\be)^{\sharp})(p^*\al)^{\sharp},(p^*\be)^{\sharp}\big)$
makes sense. 
Moreover, all operations used in the proof of \thetag{\ref{nmb:1.2}} work again, so this result holds. 
The proof in \thetag{\ref{nmb:1.6}} works and we can conclude the following result:

\begin{theorem*}
Let $p:(E,g_E)\to (B,g_B)$ be a Riemann submersion between infinite
dimensional robust Riemann manifolds. Then for 1-forms
$\al,\be\in\Om_{g_B}^1(B)$ O'Neill's formula holds in the form:
$$\boxed{
\begin{aligned}
& g_B\big(R^B(\al^{\sharp},\be^{\sharp})\be^{\sharp},\al^{\sharp}\big)
\\&\qquad
= g_E\big(R^E((p^*\al)^{\sharp},(p^*\be)^{\sharp})(p^*\be)^{\sharp},(p^*\al)^{\sharp}\big)
+\tfrac34\|[(p^*\al)^{\sharp},(p^*\be)^{\sharp}]^{\text{ver}}\|_{g_E}^2
\end{aligned}
}$$
\end{theorem*}

\section{\label{nmb:3} The diffeomorphism group $\on{Diff}_{\mathcal S}(N)$}

\subsection{\label{nmb:3.1} Diffeomorphism groups}
Let $N$ be one of the following:
\begin{itemize}
\item $N$ is a compact manifold: Then let $\on{Diff}(N)$ be the regular Lie group 
       \cite[section~38]{KM97} consisting of all smooth diffeomorphisms of $M$. 
\item $N$ is $\mathbb R^n$: we let $\on{Diff}_{\mathcal S}(\mathbb R^n)$ denote the group of all 
       diffeomorphisms of $\mathbb R^n$ which decay rapidly towards the identity. This is a regular 
       Lie group (for $n=1$ this is proved in \cite[6.4]{MGeomEvol}; the proof there works for 
       arbitrary $n$). Its Lie algebra is the space $\X_{\mathcal S}(\mathbb R^n)$ of rapidly 
       falling vector fields, with the negative of the usual bracket as Lie bracket.     
\item More generally, $(N,g)$ is a non-compact Riemannian manifold of bounded geometry, 
       see \cite{Eichhorn2007}. It is a complete Riemannian manifold and all covariant derivatives 
       of the curvature are bounded with respect to $g$. Then there is a well developed theory of 
       Sobolev spaces on $N$; let $H^\infty$ denote the intersection of all Sobolev spaces which 
       consists of smooth functions (or sections). Even on $N=\mathbb R$ the space $H^\infty$ is 
       strictly larger than the subspace $\mathcal S$ of all rapidly decreasing functions (or 
       sections) which can be defined by the condition that the Riemannian norm of all iterated 
       covariant derivatives decreases faster than the inverse of any power of the Riemannian 
       distance. There is nearly no information available on the space $\mathcal S$ for a general 
       Riemannian manifold of bounded geometry. 
			 For the following we let $\mathcal S$ denote either $H^\infty$ or the space of rapidly 
       decreasing functions. We let $\on{Diff}_{\mathcal S}(N)$ denote the group of all 
       diffeomorphisms which decay rapidly towards the identity (or differ from the identity by 
       $H^\infty$). It is a regular Lie group with Lie algebra the space $\X_{\mathcal S}(N)$ of 
       rapidly decreasing vector fields with the negative of the usual bracket. In 
       \cite[6.4]{MGeomEvol} this was proved for $N=\mathbb R$, but the same proof works for the 
       general case discussed here.     
\end{itemize}
In general, we need to impose some boundary conditions near infinity for groups of diffeomorphisms 
on a non-compact manifold $M$: The full group $\on{Diff}(N)$ of all diffeomorphisms with its natural compact 
$C^\infty$ topology is not locally contractible, so it does not admit any atlas of open charts. 

For uniformity of notation, we shall denote by $\on{Diff}_{\mathcal S}(N)$ 
any of these regular Lie groups. 
Its Lie algebra is denoted by $\X_{\mathcal S}(N)$ in each of these cases, 
with the negative of the usual bracket as Lie bracket. 
We also shall denote by $\mathcal O= C^\infty\cap \mathcal S'$ 
the space of smooth functions in the dual space $\mathcal S'$ 
(to be specific, this is the space $\mathcal O_M$ in the sense of Laurent Schwartz, 
if $N=\mathbb R^n$). 

\subsection{\label{nmb:3.2} Riemann metrics on the diffeomorphism group}
Motivated by the concept of robust Riemannian manifolds and by \cite[chapter 12]{Younes10} we will construct a right invariant weak Riemannian metric by assuming that we have a Hilbert space $\mathcal H$ together with two bounded injective linear mappings
\begin{equation*}
\X_S(N) = \Ga_{\mathcal S}(TN) \East{j_1}{} 
\mathcal H \East{j_2}{} \Ga_{C^2_b}(TN) 
\tag{1}\end{equation*}
where $\Ga_{C^2_b}(TN)$ is the Banach space of all $C^2$ vector fields $X$ on $N$ which are globally bounded together with $\nabla^gX$ and $\nabla^g\nabla^g X$ with respect to $g$, such that 
$j_2\o j_1:\Ga_{\mathcal S}(TN) \to \Ga_{C^2_b}(TN)$ is the canonical embedding. We also assume that $j_1$ has dense image.

Dualizing the Banach spaces in equation (1) and using the canonical isomorphisms between $\mathcal H$ and its dual ${\mathcal H}'$ -- which we call $L$ and $K$, we get the diagram:
\begin{equation*}\xymatrix{
\Ga_{\mathcal S}(TN) \ar@{^{(}->}[d]^{j_1} & \Ga_{{\mathcal S}'}(T^*N) \\
{\mathcal H} \ar@{^{(}->}[d]^{j_2} \ar@<.5ex>[r]^L &{\mathcal H}' \ar@{^{(}->}[u]^{j_1'} \ar@<.5ex>[l]^K\\
\Ga_{C^2_b}(TN) & \Ga_{M^2}(T^*N) \ar@{^{(}->}[u]^{j_2'}}\tag{2}\end{equation*}
Here we have written $\Ga_{{\mathcal S}'}(T^*N)$ for the dual of the space of smooth vector fields $\Ga_{\mathcal S}(TN) = \X_{\mathcal S}(N)$. We call these {\it 1-co-currents} as 1-currents are elements in the dual of $\Ga_{\mathcal S}(T^*N)$. It contains smooth measure valued cotangent vectors on $N$ (which we will write as $\Ga_{\mathcal S}(T^*N\otimes \on{vol}(N))$) and as well as the bigger subspace of second derivatives of finite measure valued 1-forms on $N$ which we have written as $\Ga_{M^2}(T^*N)$ and which is part of the dual of $\Ga_{C^2_b}(TN)$. In what follows, we will have many momentum variables with values in these spaces.

The restriction of $L$ to $\X_S(N)$ via $j_1$ gives us a positive definite weak inner product on $\X_S(N)$ which may be defined by a distribution valued kernel -- which we also write as $L$:
\begin{align*}
&\langle \quad,\quad \rangle_L: \X_{\mathcal S}(N)\x \X_{\mathcal S}(N) \to \mathbb R, \quad \text{defined by} \\
&\langle X,Y \rangle_L = \langle j_1X,j_1Y \rangle_{\mathcal H} = \iint_{N\x N} (X(y_1)\otimes Y(y_2),L(y_1,y_2)), \\
& \qquad \text{where }L \in \Ga_{\mathcal S'}(\on{pr}_1^*(T^*N) \otimes \on{pr}_2^*(T^*N))
\end{align*}
Extending this weak inner product right invariantly over $\on{Diff}_{\mathcal S}(N)$, we get a robust weak Riemannian manifold  in the sense of \ref{nmb:1.4}.

In the case (called the {\it standard case} below) that $N=\mathbb R^n$ and that 
$$\langle X,Y \rangle_L = \int_{\mathbb R^n} 
\langle (1-A\De)^l X,Y \rangle\,dx $$
we have 
$$ L(x,y) = \Big(\frac1{(2\pi)^n}\int_{\xi\in \mathbb R^n} 
e^{i\langle \xi,x-y \rangle} (1+A|\xi|^2)^l d\xi\Big)
\sum_{i=1}^n (du^i|_x \otimes dx)\otimes(du^i|_y\otimes dy) $$
where $d\xi$,  $dx$ and $dy$ denote Lebesque measure, and where $(u^i)$ are linear coordinates on $\mathbb R^n$. Here $\mathcal H$ is the space of Sobolev $H^l$ vector fields on $N$.

Note that given an operator $L$ with appropriate properties we can reconstruct the Hilbert space $\mathcal H$ with the two bounded injective mappings $j_1, j_2$.

\noindent\emph{Construction of the reproducing kernel $K$:}
The inverse map $K$ is even nicer as it is given by a $C^2$ tensor, the reproducing kernel. To see this, note that $\Ga_{M^2}(T^*N)$ contains the measures supported at one point $x$ defined by an element $\al_x\in T_x^*N$. Then $j_2(K(j_2'(\al_x)))$ is given by a $C^2$ vector field $K_{\al_x}$ on $N$which satisfies:
\begin{equation*}
\langle K_{\al_x}, X \rangle_{\mathcal H} = \al_x(j_2X)(x)\quad\text{  for all }X\in \mathcal H, 
\al_x\in T_x^*N.
\tag{3}\end{equation*}
The map $\al_x\mapsto K_{\al_x}$ is weakly $C^2_b$, thus by \cite[theorem 12.8]{KM97} 
this mapping is strongly $\on{Lip}^1$ (i.e., differentiable and the derivative is locally Lipschitz, 
for the norm on $\mathcal H$). Since $\on{ev}_y\o K: T_x^*N\ni \al_x\mapsto K_{\al_x}(y)\in T_yN$ 
is linear we get a corresponding element $K(x,y)\in L(T_x^*N,T_yN)=T_xN\otimes T_yN$ with $K(y,x)(\al_x) = K_{\al_x}(y)$. 

Using \thetag{3} twice we have (omitting $j_2$)
$$ \be_y.K(y,x)(\al_x)=\langle K(\quad,x)(\al_x), K(\quad,y)(\be_y) \rangle_{\mathcal H} = \al_x.K(x,y)(\be_y) $$
so that: 
\begin{itemize}
        \item  $K(x,y)^\top = K(y,x):T_y^*N\to T_xN$,
        \item  $K\in \Ga_{C^2_b}(\on{pr_1}^*TN \otimes \on{pr}_2^*TN)$. 
\end{itemize}
Moreover the operator $K$ defined directly by integration
\begin{align*}
&K:\Ga_{M^2}(T^*N)\to \Ga_{C^2_b}(TN) \\
&K(\al)(y_2) = \int_{y_1\in N} (K(y_1,y_2),\al(y_1)).
\end{align*}
is the same as the inverse $K$ to $L$. In fact, by definition, they agree on sections in $\Ga_{C^2}(T^*M)$ with finite support and these are weakly dense. Hence they agree everywhere. 

We will sometimes use the abbreviations $\langle \al | K |$, $| K | \be \rangle$ and 
$\langle \al | K | \be \rangle$ for the contraction of the vector values of $K$ 
in its first and second variable against 1-forms $\al$ and $\be$. 
Often these are measure valued 1-forms so after contracting, 
there remains a measure in that variable which can be integrated.

Thus the $C^2$ tensor $K$ determines $L$ and hence $\mathcal H$ and hence the whole metric on 
$\on{Diff}_{\mathcal S}(N)$. It is tempting to start with the tensor $K$, assuming it is symmetric 
and positive definite in a suitable sense. But rather subtle conditions on $K$ are required in 
order that its inverse $L$ is defined on all infinitely differentiable vector fields. For example, 
if $N = \R$, the Gaussian kernel $K(x,y)=e^{-|x-y|^2}$ does not give such an~$L$.    

In the standard case we have 
\begin{align*}
K(x,y) &= K_l(x-y) \sum_{i=1}^n \frac\p{\p x^i} \otimes 
\frac\p{\p y^i}, \\
K_l(x) &= \frac1{(2\pi)^n}\int_{\xi\in \mathbb R^n} 
\frac{e^{i\langle \xi,x \rangle}}{(1+A|\xi|^2)^l} d\xi
\end{align*}
where $K_l$ is given by a classical Bessel function of differentiability class $C^{2l}$.

\subsection{\label{nmb:3.3} The zero compressibility limit}
Although the family of metrics above does not include the case originally studied by Arnold -- the $L^2$ metric on volume preserving diffeomorphisms -- they do include metrics which have this case as a limit. Taking $N=\R^n$ and starting with the standard Sobolev metric, we can add a divergence term with a coefficient $B$:
$$ \langle X,Y \rangle_L = \int_{\mathbb R^n} \left(\langle (1-A\De)^l X,Y \rangle + B.\text{div}(X)\text{div}(Y)\right)\,dx $$
Note that as $B$ approaches $\infty$, the geodesics will tend to lie on the cosets with respect to the subgroup of volume preserving diffeomorphisms. And when, in addition, $A$ approaches zero, we get the simple $L^2$ metric used by Arnold. This suggests that, as in the so-called `zero-viscosity limit', we should be able to construct geodesics in Arnold's metric, i.e.\ solutions of Euler's equation, as limits of geodesics for this larger family of metrics on the full group.

The resulting kernels $L$ and $K$ are no longer diagonal. To $L$, we must add
$$ B\sum_{i=1}^n \sum_{j=1}^n \Big(\frac1{(2\pi)^n}\int_{\xi\in \mathbb R^n} e^{i\langle \xi,x-y \rangle} \xi_i.\xi_j d\xi\Big)(du^i|_x \otimes dx) \otimes (du^j|_y\otimes dy). $$
It can be checked that the corresponding kernel $K$ will have the form
$$ K(x,y) = K_0(x-y)\sum_{i=1}^n  \frac\p{\p x^i} \otimes \frac\p{\p y^i} + \sum_{i=1}^n \sum_{j=1}^n (K_B)_{,ij}(x-y) \frac\p{\p x^i} \otimes \frac\p{\p y^j}$$
where $K_0$ is the kernel as above for the standard norm of order $l$ and $K_B$ is a second radially symmetric kernel on $\R^n$ depending on $B$.

\subsection{\label{nmb:3.4} The geodesic equation}
According to \cite{Arnold66}, the geodesic equation on any Lie group $G$ with a right-invariant 
metric is given as follows. Let $g(t)$ be a path in $G$ and let $u(t) = \dot{g}(t).g(t)\i = 
T(\mu^{g(t)\i})\dot g(t)$ be 
the right logarithmic derivative, a path in its Lie algebra $\mathfrak g$. Here $\mu^g:G\to G$ is 
right translation by $g$. 
Then $g(t)$ is a geodesic if and only if 
$$
\p_t u = -\ad^\top_u u.
$$
where the transposed $\ad_X^\top$ is the adjoint of $\ad_X:\mathfrak{g} \rightarrow \mathfrak{g}$ 
with respect to the metric on $\mathfrak{g}$. 

In our case the Lie algebra of $\on{Diff}_{\mathcal  S}(N)$  is the space 
$\X_{\mathcal S}(N)$ of all rapidly decreasing smooth vector fields with Lie bracket (we  write  
$\ad_XY$)  the  negative  of  the usual Lie bracket 
$\ad_XY=-[X,Y]$). Then a smooth curve $t\mapsto \ph(t)$ of diffeomorphisms 
is a geodesic for the right invariant weak Riemannian metric on $\on{Diff}_{\mathcal S}(N)$ induced 
by the weak inner product $\langle  \quad,\quad\rangle_L$ on $\X_{\mathcal S}(N)$ if and only if  
$$
\p_t u = -\ad^\top_u u.
$$
as above. Here the time dependent vector field $u$ is now given by $\p_t \ph(t) = u(t)\o \ph(t)$, 
and the transposed $\ad_X^\top$ is given by 
$$
\langle \ad_X^\top Y, Z \rangle_L = \langle Y, \ad_X Z \rangle_L = -\langle Y,[X,Z] \rangle_L.
$$
The inner product is weak; existence of $\ad_X^\top$ implies condition (1) for robustness 
of the weak Riemannian manifold $(\on{Diff}_{\mathcal S}(N), \langle\quad,\quad\rangle_L)$;
it is equivalent to the fact that the dual mapping 
$\ad_X^*:\X_{\mathcal S}(N)'\to \X_{\mathcal S}(N)'$ maps the smooth dual 
$L(\X_{\mathcal S}(N))$ to itself. We also have 
$L\o \ad_X^\top = \ad_X^* \o L$.
Using  {\it  Lie derivatives}, the computation of $\ad_X^*$ is
especially  simple.  Namely,  for  any  section  
$\om$ of $T^*N \otimes  \on{vol}$   and vector fields 
$\xi,\et \in \X_{\mathcal S}(N)$, we have:
$$ \int_N (\om, [\xi,\et]) = \int_N (\om, \L_\xi(\et)) = 
-\int_N(\L_\xi(\om),\et),$$ 
hence $\ad^*_\xi(\om) = +\L_\xi(\om)$. 
Thus the Hamiltonian version of the geodesic equation on the smooth dual 
$L(\X_{\mathcal S}(N))\subset \Ga_{C^2}(T^*N\otimes \on{vol})$ becomes
$$\boxed{ 
\p_t\al  = - \ad^*_{K(\al)}\al = - \L_{K(\al)}\al,} $$
or, keeping track of everything,
\begin{equation*}\boxed{\begin{aligned}
\p_t\ph &= u\o \ph, \\ 
\p_t\al &= - \L_u\al \\
u = K(&\al) = \al^\sharp,\quad \al=L(u) = u^\flat.
\end{aligned}
\tag{1}
}\end{equation*}
One can also derive the geodesic equation from the conserved momentum mapping 
$J:T\on{Diff}_{\mathcal S}(N)\to \X_{\mathcal S}(N)'$ given by $J(g,X) = L\o\on{Ad}(g)^\top X$ where 
$\on{Ad}(g)X = Tg\o X\o g\i$. This means that $\on{Ad}(g(t))u(t)$ is conserved and 
$0=\p_t \on{Ad}(g(t))u(t)$ leads quickly to the geodesic equation. It is remarkable that the 
momentum mapping exists if and only if $(\on{Diff}_{\mathcal S}(N), \langle  \quad,\quad\rangle_L)$ 
is a robust weak Riemannian manifold.

\section{\label{nmb:5} The differentiable Chow manifold (alias the non-linear Grassmannian)}

\subsection{\label{nmb:5.1} The differentiable Chow manifold as a homogeneous space for 
$\on{Diff}_{\mathcal S}(N)$ and the induced weak Riemannian metric} 
Let $M$ be a compact manifold with $\on{dim}(M) < \on{dim}(N)$.  
The space of submanifolds of $N$ diffeomorphic to $M$ will be called $B(M,N)$. In the case $m=0$ 
and $N=\R^D$, i.e.\ $M$ is a finite set of, say $p$, points in Euclidean $D$-space, the space 
$B(M,N)$ is what we called the space of landmark points $\L^p(\R^D)$ in our earlier paper \cite{MMM1}. 

$B(M,N)$ can be viewed as a quotient of $\on{Diff}_{\mathcal S}(N)$. If we fix a base submanifold $F_0\subset N$ 
diffeomorphic to $M$, then we get a map of $\on{Diff}_{\mathcal S}(N)$ into $B(M,N)$ by 
$\ph \mapsto \ph(F_0)$. The image will be an open subset $B_0(M,N)$ of $B(M,N)$ which is the 
quotient of $\on{Diff}_{\mathcal S}(N)$ by the subgroup of diffeomorphisms which map $F_0$ to itself.  
We will study $B(M,N)$ using this approach and without further comment 
replace the full space $B(M,N)$ by this component $B_0(M,N)$.

The normal bundle to $F \subset N$ may be defined as $TB^\bot \subset TN|_B$, with the help of an auxiliary  Riemann metric on $N$. But we want to avoid this auxiliary metric, so we shall define the normal bundle as the quotient $\on{Nor}(F):= TN|_F/TF$ over $F$. Then its dual bundle, the conormal bundle, is $\on{Nor}^*(F)=\text{Annihilator}(TF)\subset T^*N|_F$, a sub-bundle not a quotient. The tangent space $T_FB(M,N)$to $B(M,N)$ at $F$ can be identified with the space of all smooth sections $\Ga_{\mathcal S}(\on{Nor}(F))$ of the normal bundle. 

A simple way to construct local coordinates on $B(M,N)$ near a point $F \in B(M,N)$ is to 
trivialize a neighborhood of $F \subset N$. 
To be precise, assume we have a tubular neighborhood, 
i.e., an isomorphism $\Ph$:
\begin{equation*}
\begin{array}{rcl}
B(M,N)& & \on{Nor}(F) \\ 
\cup\hspace*{.1cm} & & \hspace*{.1cm}\cup \\
U_B & \hspace*{-.2cm}\stackrel{\Ph}{\longrightarrow}\hspace*{-.2cm} & U_N \\
\cup\hspace*{.1cm} & & \hspace*{.1cm}\cup \\
F\hspace*{.08cm} & \hspace*{-.2cm}=\hspace*{-.2cm} & \text{0-section} \\
\end{array}
\end{equation*} 
from an open neighborhood $U_B$ of $F$ in $N$ to an open neighborhood $U_N$ of the 0-section in the normal bundle $\on{Nor}(F)$. Assume moreover that $\Ph$ is the identity on $F$ and its normal derivative along $F$ induces the identity map on $\on{Nor}(F)$. The map $\Ph$ induces a local projection $\pi:U_B \rightarrow F$ and partial linear structure in the fibres of this projection. Then we get an open set $U_\Ph \subset B(M,N)$ consisting of submanifolds $F' \subset U_B$ which intersect the fibres of $\pi$ normally in exactly one point. Under $\Ph$ these submanifolds are all given by smooth sections of $\on{Nor}(F)$ which lie in $U_N$. If we call this set of sections $U_\Ga$ we have a chart:
$$ B(M,N) \supset U_\Ph \cong U_\Ga \subset \Ga_{\mathcal S}(\on{Nor}(F))$$

We define a Riemannian metric on $B(M,N)$ following the procedure used for $\on{Diff}_{\mathcal S}(N)$. 
For any $F \subset N$, we decompose $\mathcal H$ into:
\begin{align*}
\mathcal H^{\text{vert}}_F &= 
j_2\i\bigl(\{X\in\Ga_{C^2_b}(TN): X(x) \in T_xF, \text{ for all } x\in F\}\bigr) \\
\mathcal H^{\text{hor}}_F &= 
\text{perpendicular complement of } \mathcal H^{\text{vert}}_F
\end{align*}
It is then easy to check that we get the diagram:
\begin{equation*}\xymatrix{
\Ga_{\mathcal S}(TN) \ar@{^{(}->}^(.6){j_1}[r] \ar@{>>}[d]^{\text{res}} & \mathcal H \ar@{^{(}->}[r]^(.4){j_2} \ar@{>>}[d] & \Ga_{C^2_b}(TN) \ar@{>>}[d]^{\text{res}}\\
 \Ga_{\mathcal S}(\on{Nor}(F)) \ar@{^{(}->}[r]^(.6){j_1^f} & \mathcal H^{\text{hor}}_F \ar@{^{(}->}[r]^(.4){j_2^f}  & \Ga_{C^2_b}(\on{Nor}(F)).
 }\end{equation*}
As this is an orthogonal decomposition, $L$ and $K$ take 
${\mathcal H}^{\text{vert}}_F$ and ${\mathcal H}^{\text{hor}}_F$ 
into their own duals and, as before we get:
\begin{equation*}\xymatrix{
\Ga_{\mathcal S}(\on{Nor}(F)) \ar@{^{(}->}[d]^{j_1} & \Ga_{{\mathcal S}'}(\on{Nor}^*(F)) \\
{\mathcal H}^{\text{hor}}_F \ar@{^{(}->}[d]^{j_2} \ar@<.5ex>[r]^{L_F} &({\mathcal H}^{\text{hor}}_F)' \ar@{^{(}->}[u]^{j_1'} \ar@<.5ex>[l]^{K_F} \\
\Ga_{C^2_b}(\on{Nor}(F)) & \Ga_{M^2}(\on{Nor}^*(F)) \ar@{^{(}->}[u]^{j_2'} }\end{equation*}
$K_F$ is just the restriction of $K$ to this subspace of ${\mathcal H}'$ and is given by the kernel:
$$K_F(x_1,x_2):= \text{image of }K(x_1,x_2) \in \on{Nor}_{x_1}(F)\otimes \on{Nor}_{x_2}(F)), \quad x_1, x_2 \in F.$$ 
This is a $C^2$ section over $F\times F$ of $\on{pr}_1^*\on{Nor}(F)
\otimes \on{pr}_2^*\on{Nor}(F)$. We can identify the space of horizontal vector fields ${\mathcal H}^{\text{hor}}_F$ as the closure of the image under $K_F$ of measure valued 1-forms supported by $F$ and with values in $\on{Nor}^*(F)$.
A dense set of elements in ${\mathcal H}^{\text{hor}}_F$ is given by either taking the 1-forms with finite support or taking smooth 1-forms. In the first approach, ${\mathcal H}^{\text{hor}}_F$ is the closure of the span of the vector fields 
$\big| K_F(\cdot,x) \big | \al_x \big\rangle$ where $x \in F$ and $\al_x \in \on{Nor}_x^*(F)$. In the smooth case, fix a volume form $\ka$ on $M$ and a smooth covector 
$\xi \in \Ga_{\mathcal S}(\on{Nor}^*(F))$.   Then $\xi.\ka$ defines a horizontal vector field $h$ like this:
$$h(x_1) =  \int_{x_2\in F} \big| K_F(x_1,x_2)\big| \xi(x_2).\ka(x_2) \big\rangle$$
The horizontal lift $h\h$ of any $h \in T_FB(M,N)$ is then: 
$$ h\h (y_1) = K(L_F h)(y_1) = 
\int_{x_2\in F}\big| K(y_1,x_2)\big| L_F h(x_2)\big\rangle, \quad y_1 \in N.$$
Note that all elements of the cotangent space $\al \in \Ga_{\mathcal S'}(\on{Nor}^*(F))$ can be pushed up to $N$ by $(j_F)_*$, where $j_F:F \hookrightarrow N$ is the inclusion, and this identifies $(j_F)_*\al$ with a 1-co-current on $N$.

Finally the induced homogeneous weak Riemannian metric on $B(M,N)$ is given like this:
\begin{align*}
\langle  h,k\rangle_F &= \int_N (h\h (y_1),L(k\h )(y_1)) 
= \int_{y_1\in N} (K(L_F h))(y_1),(L_F k)(y_1)) \\
& = \int_{(y_1,y_2)\in N\x N} (K(y_1,y_2),(L_F h)(y_1)\otimes 
(L_F k)(y_2)) \\
& = \int_{(x_1,x_1)\in F\x F} \big\langle L_F h(x_1) \big| 
K_F(x_1,x_2) \big| L_F h(x_2) \big\rangle  
\end{align*}
With this metric, the projection from $\on{Diff}_{\mathcal S}(N)$ to $B(M,N)$ is a submersion. 
The inverse co-metric on the smooth cotangent bundle 
$\bigsqcup_{F\in B(M,N)}\Ga(\on{Nor}^*(F)\otimes \on{vol}(F))\subset T^*B(M,N)$
is much simpler and easier to handle:
$$\boxed{
\langle \al,\be \rangle_F = \iint_{F\x F} \big\langle \al(x_1) \big| K_F(x_1,x_2)\big|\be(x_1)\big\rangle.
}$$
It is simply the restriction to the co-metric on the 
Hilbert sub-bundle of $T^*\on{Diff}_{\mathcal S}(N)$ defined by $\mathcal H'$ 
to the Hilbert sub-bundle of subspace $T^*B(M,N)$ defined by $\mathcal H_F'$.

Because they are related by a submersion, the geodesics on $B(M,N)$ are the horizontal geodesics  
on $\on{Diff}_{\mathcal S}(N)$, as described in box (1), section 3.4. We have two variables: a 
family $\{F(t)\}$ of submanifolds in $B(M,N)$ and a time varying momentum $\al(t,\cdot) 
\in \on{Nor}^*(F(t))\otimes \on{vol}(F(t))$ which lifts to the horizontal 1-co-current 
$(j_{F(t)})_*(\al(t,\cdot)$ on $N$. Then the horizontal geodesic on $\on{Diff}_{\mathcal S}(N)$ is given by the 
same equations as before: 
\begin{equation*}\boxed{\begin{aligned}
\p_t (F(t)) &= \on{res}_{\on{Nor}(F(t))}(u(t,\cdot))\\
u(t,x) &= \int_{(F(t))_y} \big| K(x,y)\big| \al(t,y) \big\rangle \in \X_{\mathcal S}(N)\\
\p_t\left((j_{F(t)})_*(\al(t,\cdot)\right) &= -{\mathcal L}_{u(t,\cdot)}((j_{F(t)})_*(\al(t,\cdot)).
\end{aligned}} \end{equation*}
This is a complete description for geodesics on $B(M,N)$ but it is not very clear how to compute the Lie derivative of $(j_{F(t)})_*(\al(t,\cdot)$. One can unwind this Lie derivative via a torsion-free connection, but we turn to a different approach which will be essential for working out the curvature of $B(M,N)$.

\subsection{\label{nmb:5.2} Auxiliary tensors on $B(M,N)$}
Our goal is to reduce calculations on the infinite dimensional space $B(M,N)$ to calculations on the finite dimensional space $N$. To do this we need to construct a number of useful tensors on $B(M,N)$ from tensors on $N$ and compute the standard operations on them. These will enable us to get control of the geometry of $B(M,N)$. Let $m$ be the dimension of $M$, $n$ the dimension of $N$. For $F \in B(M,N)$, let $j_F:F \hookrightarrow N$ be the embedding. We will assume that $M$ is orientable for simplicity, so that $\on{vol}(M) \cong \Om^m(M)$.

\noindent\thetag{\nmb:{1}}
We denote by $\ell$ the left action:
$$ \ell: \on{Diff}_{\mathcal S}(N) \times B(M,N) \rightarrow B(M,N)$$
given by $\ell(\ph,F) \text{ or } \ell^F(\ph)= \ph(F)$. 
For a vector field $X\in\X_{\mathcal S}(N)$ let $B_X$ be the infinitesimal action 
(or fundamental vector field) on $B(M,N)$ given by $B_X(F)= T_{\on{Id}}(\ell^F)X$ with its flow 
$\Fl^{B_X}_t(F) = \Fl^X_t(F)$. 
The fundamental vector field mapping of a left action is a 
Lie algebra anti-homomorphism and the Lie bracket on $\on{Diff}_{\mathcal S}(N)$ 
is the negative of the usual Lie bracket on $\X_{\mathcal S}(N)$, 
so we have $[B_X,B_Y]=B_{[X,Y]}$. 
The set of these vectors $\{B_X(F): X\in \X_{\mathcal S}(N)\}$ equals the whole tangent space $T_FB(M,N)$.

\noindent\thetag{\nmb:{2}}
Note that $B(M,N)$ is naturally submanifold of the vector space of $m$-currents on $N$:
$$ B(M,N) \hookrightarrow \Om^m_{\mathcal S}(N)'=\Ga_{\mathcal S'}(\La^mTN), \quad \text{via } F 
\mapsto \left(\om \mapsto \int_F\om\right).$$
Any $\al\in\Om^m(N)$ is a linear coordinate on $\Ga_{\mathcal S'}(\La^mTN)$ 
and this restricts to the function $B_\al\in C^\infty(B(M,N),\mathbb R)$ 
given by $B_\al(F)= \int_F\al$. If $\al=d\be$ for $\be\in\Om^{m-1}(N)$ then 
$$B_\al(F)=B_{d\be}(F)=\int_Fj_F^*d\be= \int_Fdj_F^*\be = 0$$ 
by Stokes' theorem. 

For $\al\in\Om^m(N)$ and $X\in\X_{\mathcal S}(N)$ we can evaluate the vector 
field $B_X$ on the function $B_\al$:
\begin{align*}
B_X(B_\al)(F) &= dB_\al(B_X)(F) = \p_t|_0 B_\al(Fl^X_t(F))= 
\int_F j_F^*\L_X\al = B_{\L_X(\al)}(F) \\
& \text{as well as } = \int_F j_F^*(i_X d\al+d i_X\al) = 
\int_F j_F^* i_X d\al = B_{i_X(d\al)}(F)
\end{align*}
If $X\in\X_{\mathcal S}(N)$ is tangent to $F$ along $F$ then
$B_X(B_\al)(F) = \int_F \L_{X|_F}j_F^*\al = 0$.

More generally, a $pm$-form $\al$ on $N^k$ defines a function $B^{(p)}_\al$ on $B(M,N)$ 
by $B^{(p)}_\al(F) = \int_{F^p} \al$. Using this for $p=2$, 
we find that for any two $m$-forms $\al, \be$ on $N$, the inner product of $B_\al$ and $B_\be$ is given by:
$$ g\i_B(B_\al, B_\be) = B^{(2)}_{\langle \al | K | \be \rangle}.$$

\noindent\thetag{\nmb:{3}}
For $\al\in\Om^{m+k}(N)$ we denote by $B_\al$ the $k$-form in $\Om^k(B(M,N))$ 
given by the skew-symmetric multi-linear form: 
$$(B_\al)_F(B_{X_1}(F),\dots,B_{X_k}(F))= \int_F j_F{}^*( i_{X_1\wedge \dots\wedge X_k}\al).$$
This is well defined: If one of the $X_i$ is tangential to $F$ at a point $x\in F$ then $j_F{}^*$ 
pulls back the resulting $m$-form to 0 at $x$. 

Note that any smooth cotangent vector $a$ to $F \in B(M,N)$ is equal to 
$B_\al(F)$ for some closed $(m+1)$-form $\al$. 
Smooth cotangent vectors at $F$ are elements of 
$\Ga_{\mathcal S}(F,\on{Nor}^*(F) \otimes \Om^m(F))$. 
Fix a nowhere zero global section $\ka$ of $\Om^m(F)$. 
Then $\frac{a}{\ka}$ is the differential of a unique function $f$ 
on the normal bundle to $F$ which is linear on each fibre. 
Let $\ph$ be a local isomorphism from a neighborhood of $F$ in $N$ 
to a neighborhood of the 0-section in this normal bundle and let 
$\rh$ be a function on the normal bundle which is one near the 0-section and 
has support in this neighborhood. 
Take $\al = d(f.\ka \circ \ph)$ (extended by zero). It's easy to see that this does it. 

Likewise, a $pm+k$ form 
$\al \in \Om^{pm+k}(N^p)$ defines a $k$-form on $B(M,N)$ as follows:
First, for $X\in\X_{\mathcal S}(N)$ let $X^{(p)}\in\X(N^p)$ be given by  
\begin{multline*}
X^{(p)}_{(n_1,\dots,n_p)} := (X_{n_1}\x 0_{n_2}\x\dots\x 0_{n_p})  + 
(0_{n_1}\x X_{n_2}\x 0_{n_3}\x\dots\x 0_{n_p}) + \dots
\\ \dots + (0_{n_1}\x\dots\x0_{n_{p-1}}\x X_{n_p}).
\end{multline*}
Then we put
$$(B_\al^{(p)})_F(B_{X_1}(F),\dots,B_{X_k}(F))
= \int_{F^p} j_{F^p}{}^*( i_{X_1^{(p)}\wedge \dots\wedge X_k^{(p)}}\al).$$
This is just $B$ applied to the submanifold $F^p\subset N^p$ and to the special vector fields
$X^{(p)}$. Thus all properties of $B$ continue to hold for $B^{(p)}$; in particular, \thetag{4} 
below hold for $X^{(p)}$ instead of $X$.

\noindent\thetag{\nmb:{4}}
We have $\boxed{i_{B_X}B_\al= B_{i_X\al}}$ because
\begin{align*}
\big(i_{B_{X_1}}B_\al\big)&\big(B_{X_2},\dots,B_{X_k}\big)(F)
=B_\al\big(B_{X_1},B_{X_2},\dots,B_{X_k}\big)(f)
\\&
=\int_F j_F{}^*\big(i_{X_k}\dots i_{X_2} (i_{X_1}\al)\big)
=B_{i_{X_1}\al}\big(B_{X_2},\dots,B_{X_k}\big)(F)
\end{align*}
For the exterior derivative we have
$\boxed{\;dB_\al = B_{d\al}\;}$ for any $\al\in\Om^{m+k}(N)$. Namely,
\begin{align*}
(&dB_\al)(B_{X_0},\cdot\cdot,B_{X_k})(F)
= \sum_{i=0}^k (-1)^i B_{X_i}(B_\al(B_{X_0},\cdot\cdot,\widehat{B_{X_i}},\cdot\cdot,B_{X_k}))(F)
\\&\qquad
+ \sum_{i<j} (-1)^{i+j} B_\al(B_{[X_i,X_j]},B_{X_0},\cdot\cdot,\widehat{B_{X_i}},\cdot\cdot\widehat{B_{X_j}},\cdot\cdot,B_{X_k}))(F)
\\&
= \sum_{i=0}^k (-1)^i \int_F j_F^*i_{X_i} d i_{X_0\wedge \cdot\cdot\widehat{X_i}\cdot\cdot\wedge X_k} \al
+ \sum_{i<j} (-1)^{i+j} \int_F j_F^* 
i_{[X_i,X_j]\wedge X_0\wedge \cdot\cdot\widehat{X_i}\cdot\cdot\widehat{X_j}\cdot\cdot\wedge X_k}\al
\\&
= \int_F j_F{}^*\Big(\sum_{i=0}^k (-1)^i \L_{X_i}  i_{X_k} \cdot\cdot\widehat{i_{X_i}}\cdot\cdot i_{X_0} 
\\&\qquad\qquad\quad
- \sum_{i<j} (-1)^{i} 
i_{X_0\wedge \cdot\cdot\widehat{X_i}\cdot\cdot\wedge{X_{j-1}}\wedge [X_i,X_j]\wedge X_{j+1}\cdot\cdot\wedge X_k}\Big)\al
\\&
= \int_F j_F{}^*\sum_{i=0}^k (-1)^i \Big(\L_{X_i}  i_{X_k} \hspace*{-.1cm}\cdot\cdot\widehat{i_{X_i}}\hspace*{-.1cm}\cdot\cdot i_{X_0} 
\\&\qquad\qquad\qquad\qquad\qquad
- \sum_{j=i+1}^k  
i_{X_k}\hspace*{-.1cm}\cdot\cdot i_{X_{j+1}}\,[\L_{X_i},i_{X_j}]\,i_{X_{j-1}}\hspace*{-.1cm}\cdot\cdot\widehat{i_{X_i}}\hspace*{-.1cm}\cdot\cdot i_{X_0}\Big)\al
\\&
= \int_F j_F{}^*\Big(\sum_{i=0}^k (-1)^i  i_{X_k} \hspace*{-.1cm}\cdot\cdot i_{X_{i+1}} \L_{X_i} i_{X_{i-1}}\hspace*{-.1cm}\cdot\cdot i_{X_0}\al\Big) 
\\&
= \int_F j_F{}^*\Big(\sum_{i=0}^k (-1)^i  i_{X_k} \hspace*{-.1cm}\cdot\cdot i_{X_{i+1}}(d\,i_{X_i}+i_{X_i}d) i_{X_{i-1}}\hspace*{-.1cm}\cdot\cdot i_{X_0}\al\Big) 
\\&
= \int_F j_F{}^*\Big(\sum_{i=0}^k (-1)^i  i_{X_k} \cdot\cdot i_{X_{i+1}}d\,i_{X_i}\cdot\cdot i_{X_0}
+\sum_{i=0}^k (-1)^i  i_{X_k} \hspace*{-.1cm}\cdot\cdot i_{X_i}d\, i_{X_{i-1}}\hspace*{-.1cm}\cdot\cdot i_{X_0}\Big)\al 
\\&
= 0+ \int_F j_F{}^* i_{X_k} \cdot\cdot i_{X_0}d\al = B_{d\al}(B_{X_0},\cdot\cdot,B_{X_k})(F)
\end{align*}
Finally we have
$\boxed{
\L_{B_X}B_\al = B_{\L_X\al}
}$
since
\begin{align*}
\L_{B_X}B_\al &= (i_{B_X}\,d + d\,i_{B_X})B_{\al}
=B_{(i_Xd+di_X)\al} = B_{\L_X\al}.
\end{align*}
Note that these identities generalize the results in item (2).

\noindent\thetag{\nmb:{5}}
For $\al\in\Om^{m+1}(N)$ we pull back to $\on{Diff}_{\mathcal S}(N)$ the 1-form $B_{\al}$ on $B(M,N)$ 
 where $\ph_0(F_0)= F$:
$$\boxed{
\begin{aligned}
\bigl((\ell^{F_0})^*B_{\al}\bigr)_{\ph_0} (X\o\ph_0) &= \bigl((\ell^F)^*B_\al\bigr)_{\on{Id}}(X) 
= (B_{\al})_F(B_X(F)) = \int_F j_F^*i_X\,\al,
\\
\bigl((\ell^{F})^*B_{\al}\bigr)_{\on{Id}} &= \al|_F =: \mu(\al,F) = \mu_\al(F) = \mu^F(\al) \in \X_{\mathcal S}(N)'
\\ \mu:\Om^{m+1}(N) &\x B(M,N) \to \X_{\mathcal S}(N)'
\\
\mu(\al,F)\quad&\text{is a 
1-cocurrent with support along }F.
\end{aligned}
}$$
The mapping $\mu:\Om^{m+1}(N) \x B(M,N) \to \X_c(N)'$ is smooth, $\mu^F:\Om^{m+1}(N)\to \X_c(N)'$ is 
bounded linear, and the differential of $\mu_\al:B(M,N)\to\X_{\mathcal S}(N)'$
is computed as follows:
\begin{align*}
\langle d(\mu_\al)&(B_X(F)),Y\rangle 
=\langle D_{F,B_X}\mu(\al,F),Y\rangle 
= D_{F,B_X} \langle \mu(\al,F),Y\rangle 
= \p_t|_0 \langle \al_{\Fl^X_t(F)},Y \rangle 
\\&
= \p_t|_0 \int_{\Fl^X_t(F)} j_{\Fl^X_t(F)}{}^*i_Y\al
= \p_t|_0 \int_{\Fl^X_t(F)} (\Fl^X_t\o j_F\o (\Fl^X_t|_F)\i)^*i_Y\al
\\&
= \p_t|_0 \int_{\Fl^X_t(F)} (\Fl^X_t|_F)\i)^*j_F{}^*(\Fl^X_t)^* i_Y\al
\\&
= \p_t|_0 \int_{F} j_F{}^*(\Fl^X_t)^* i_Y\al
=  \int_{F} j_F{}^*\,\L_X(i_Y\al)
\\&
=  \int_{F} j_F{}^*(i_{[X,Y]}\al + i_Y\L_X\al)
=  \langle \mu(\al,F),\L_XY \rangle + \langle \mu(\L_X\al,F), Y \rangle.
\end{align*}
This means
\begin{equation*}\boxed{
d\mu_\al(B_X(F)) = \mu(\al,F)\o \L_X  +\mu(\L_X\al,F) = -\L_X\mu(\al,F) + \mu(\L_X\al,F),
}\tag{\nmb:{6}}\end{equation*}
where $\L_X\mu(\al,F)$ denotes the Lie derivative of 1-currents.
There are two interpretations of formula \thetag{\nmb|{6}}:
\begin{align*}
d\mu_\al(B_X) &= -\L_X\o \mu_\al + \mu_{\L_X\al},
\\
d\mu_\al(B_X(F)) &= -(\L_X \mu^F)(\al). 
\end{align*}
We shall also need the mapping
$\mu:\Om^m(N)\x B(M,N) \to C^\infty_c(N)'$ with values in the linear space of 
distributions (without the density part) on $N$ which is given by
$$
\langle \mu(\ga,F),f \rangle = \int_F f.\ga = \int_F j_F{}^* (g\ga).
$$
The distribution $\mu(\ga,F)$ is again bounded linear in $\ga\in\Om^m(N)$, and its derivative with 
respect to $F$ is given by \thetag{6} again, with the same proof as above. 

\section{\label{nmb:5}Geodesics and curvature on $B(M,N)$}

We want to use the auxiliary tensors of the last section to derive formulas for geodesics and 
curvature on $B(M)$, using Mario's formula to compute the curvature. The basic idea is to write a 
smooth co-vector $a$ at a point $F \in B(M,N)$ as $B_\al$ where $\al$ is an $(m+1)$-form on $N$. As 
always, for any $(m+1)$-form $\al$ on $N$, $B_\al^\sh$ is the ($C^2$) vector field on $B(M,N)$ which 
is dual to the smooth 1-form $B_\al$. At each point $F \in B$, $B_\al^\sh$ lifts horizontally to a 
tangent vector at the identity to $\on{Diff}_{\mathcal S}(N)$, which is given by the vector field
$$\mu(\al,F)^\sh = \int_N |K|\mu(\al,F)\rangle \in\X_{C^2}(N)$$ 
so that $B_{\mu(\al,F)^\sh}(F)=B_\al^\sh(F)$. See (\ref{nmb:5.2}.5).

With these co-vectors, we consider next the force introduced in section 2.3. We have:
$$ 2\mathcal F(\al,\be)=d(\langle B_\al, B_\be \rangle) 
= d\left(B^{(2)}_{\langle \al | K | \be \rangle}\right) 
= B_{d(\langle \al | K | \be \rangle)}^{(2)}.$$
But $\langle \al | K | \be \rangle$ is a $2m$-form on $N \times N$ 
and $d$ can be split into two parts $d_1+d_2$ acting on the first and second factors. Evaluating 
this 1-form at $F$ and taking its inner product with $B_X, X \in \X_{\mathcal S}(N)$, we get: 
\begin{align*}
\Big(B_{d(\langle \al | K | \be \rangle)}^{(2)}&(F),B_X(F)\Big) =
\iint_{F\times F} j_{F\x F}{}^* i_{X^{(2)}} (d(\langle \al | K | \be \rangle)) \\
&= \iint_{F\times F} j_{F\x F}{}^* \big((i_X)_1 (d_1(\langle \al | K | \be \rangle))
+(i_X)_2 (d_2(\langle \al | K | \be \rangle))\big) \\
\text{because } &F\times F \text{ has type } (m,m) \text{ and the integrand must have the same type}\\
&= \int_F j_F^* i_X d\left(i_{\mu(\be,F)^\sh}(\al) + i_{\mu(\al,F)^\sh}(\be)\right)
\end{align*}
hence  
$$ 2\mathcal F(\al,\be) = B_{d(\langle \al | K | \be \rangle)}^{(2)} 
= B_\ga, \quad \ga = \L_{\mu(\be,F)^\sh}(\al) + \L_{\mu(\al,F)^\sh}(\be).$$
Here the superscript 2 on $X$ means that $X^{(2)}$ is the vector field on $N \times N$ given by $0 \times X + X \times 0$ whereas on $B$, because $d(\langle \al | K | \be \rangle)$ is a $(2m+1)$-form on $N \times N$, we must apply $B^{(2)}$, not $B$, to it. Thus we define the force $F$ using operations on the finite dimensional manifold $N$ by:
$$\boxed{ \mathcal F_N(\al,\be,F) 
:= \big(\text{image in Nor$^*(F)\otimes$vol}(F)\big)\left(\frac{1}{2}(\L_{\mu(\be,F)^\sh}(\al) 
+ \L_{\mu(\al,F)^\sh}(\be))\right).}$$
The term `force' comes from the fact that the geodesic acceleration is given by $\mathcal F(\al,\al)$. In our case, we find that the geodesic equation on $B(M,N)$ can be extended to an equation in the variables $F(t) \in B(M,N)$ 
and $\al(t,\cdot)$ a time varying $(m+1)$-form on $N$:
$$\boxed{\begin{aligned}
\p_t (F(t)) &= \text{(res to Nor}(F))u\\
u &= \mu(\al,F)^\sh = \int_{F(t)(y)} | K(\cdot,y) | \al(y) \rangle \\
\p_t(\al) &= \mathcal F(\al,\al,F) = \L_u(\al).
\end{aligned}}$$

Moving to curvature, fix $F$. Then we claim that for any two smooth co-vectors $a, b$ at $F$, we can construct not only two closed $(m+1)$-forms $\al, \be$ on $N$ as above but also two commuting vector fields $X_\al, X_\be$ on $N$ in a neighborhood of $F$ such that:
\begin{enumerate}
\item $B_\al(F)=a$ and $B_\be(F)=b$,
\item $B_{X_\al}(F)=a^\sh$ and $B_{X_\be}(F)=b^\sh$
\item $\L_{X_\al}(\al) = \L_{X_\al}(\be) = \L_{X_\be}(\al) = \L_{X_\be}(\be) = 0$
\item $[X_\al, X_\be] = 0$
\end{enumerate}
We can do this using a local isomorphism of $N$ with the normal bundle 
to $F$ in $N$ as above. This gives a projection $\pi$ of a neighborhood of $F$ in $N$ to $F$ and partial linear structure on its fibres. Then for $\al$ and $\be$ use $(m+1)$-forms $\ka \wedge \om$ where $\ka$ is a pull back of an $m$-form on $F$ and $\om$ is a 1-form constant along the fibres; and for $X_\al$ and $X_\be$ use vector fields which are tangent to the fibres of $\pi$ and constant with respect to the linear structure on them. 

We are now in a position to use the version (\ref{eq:curv3}) of Mario's formula. As it stands, this formula calculates curvature using operations on $B(M,N)$. What we want to do is to write everything using forms and fields on $N$ instead. We first need an expression for the stress $\mathcal D(\al,\be)$ in this formula. Using notation from (\ref{nmb:1.3}.2):
\begin{align*}
\mathcal D(\al,\be,F) &= D_{F,B_{X_\al}(F)}(B^\sh_\be - B_{X_\be})\\
&= [B_{X_\al},B^\sh_\be - B_{X_\be}](F) = [B_{X_\al},B^\sh_\be](F).
\end{align*}
In order to compute the Lie bracket, we apply it to a smooth function 
$B_\ga$ on $B(M,N)$ where $\ga\in\Om^m(N)$. 
Then we have, using \ref{nmb:5.2} repeatedly:
\begin{align*}
&(\L_{B_\be^\sh}B_\ga)(F) = (\L_{B_{\mu(be,F)^\sh}}B_\ga)(F) = B_{\L_{\mu(be,F)^\sh}\ga}(F)
\\&
(\L_{B_{X_\al}}\L_{B_\be^\sh}B_\ga)(F) = (\L_{B_{X_\al}} B_{\L_{\mu(\be,F)^\sh}\ga})(F)
\\& \qquad
= B({\L_{D_{F,B_{X_\al}}\mu(\be,F)^\sh}\ga})(F) +  B(\L_{X_\al} \L_{\mu(\be,F)^\sh}\ga)(F)
\\&
(\L_{B_\be^\sh}\L_{B_{X_\al}}B_\ga)(F) = (\L_{B_{\mu(\be,F)^\sh}}B_{\L_{X_\al}\ga})(F) 
= B(\L_{\mu(\be,F)^\sh}\L_{X_\al}\ga)(F)
\\&
D_{F,B_{X_\al}} \mu(\be,F)\s  
= D_{F,B_{X_\al}} \int_N |K|\mu(\be,F)\rangle 
= \int_N|K| D_{F,B_{X_\al}} \mu(\be,F)\rangle
\\&\qquad
=\int_N\big|K\big|(-\L_{X_\al}\mu(\be,F) + \mu(\L_{X_\al}\be,F))\big\rangle
\quad\text{  by \thetag{\ref{nmb:5.2}.6}}
\\&\qquad
=\int_N\big|\L_{0\x X_\al}K\big|\mu(\be,F)\big\rangle 
+ \mu(\L_{X_\al}\be,F)\s
\\&
([B_{X_\al},B_{\mu(be,F)^\sh}]B_\ga)(F) =  (\L_{B_\be^\sh}B_\ga - \L_{B_\be^\sh}\L_{B_{X_\al}}B_\ga)(F) 
\\&\qquad
= B({\L_{D_{F,B_{X_\al}}\mu(\be,F)^\sh}\ga})(F) +  B(\L_{[X_\al,\mu(\be,F)^\sh]}\ga)(F)
\\&\qquad
= (\L_{B(D_{F,B_{X_\al}}\mu(\be,F)^\sh+[X_\al,\mu(\be,F)^\sh])} B_\ga)(F)
\\&
[B_{X_\al},B_{\mu(be,F)^\sh}](F) = B(D_{F,B_{X_\al}}\mu(\be,F)^\sh+\L_{X_\al}\mu(\be,F)^\sh)
\\&
= B\Big(\int_N\big|\L_{0\x X_\al}K\big|\mu(\be,F)\big\rangle 
+ \mu(\L_{X_\al}\be,F)\s+
+\int_N\big|\L_{X_\al\x 0}K\big|\mu(\be,F)\big\rangle\Big)
\\&\qquad
= B\Big(\int_N\big|\L_{X_\al{}^{(2)}}K\big|\mu(\be,F)\big\rangle\Big) + 0. 
\end{align*}

Thus we define the stress $\mathcal D=\mathcal D_N$ on $N$ by:
$$\boxed{ \mathcal D(\al,\be,F)(x) 
= \big(\text{restr.\ to Nor}(F)\big)
\left(-\int_{y\in F} \Big|\L_{X_\al^{(2)}}(x,y) K(x,y) \Big| \be(y) \Big\rangle\right).}$$

Next consider the second derivative terms in $R_{11}$. A typical term works out as follows:
\begin{align*}
B_{X_\al}B_{X_\al}(<B_\be, B_\be>) &=\L_{B_{X_\al}}\L_{B_{X_\al}}(<B_\be, B_\be>)
= B_{\L_{X_\al^{(2)}}\L_{X_\al^{(2)}}<\be|K|\be>}
\\& 
= B_{<\be|\L_{B_{X_\al^{2}}}\L_{X_\al^{(2)}}K|\be>}
\end{align*}
Extending Lie bracket notation slightly, we can write  
$$<\be|\L_{X_\al^{(2)}}\L_{X_\al^{(2)}}K|\be>  = 
\Big\langle \be \Big|[X_\al^{(2)},[X_\al^{(2)},K]] \Big| \be \Big\rangle.$$
Analogous formulas hold for the other terms.

Finally, putting everything together, we find the formula for curvature:
$$\boxed{
\begin{aligned}
\langle R_{B(M,N)}&(B_\al^\sh, B_\be^\sh) B_\be^\sh, B_\al^\sh \rangle(F) = R_{11}+R_{12}+R_2+R_3 \\
R_{11} &= \tfrac12 \iint_{F \times F} \Big( \big\langle \be \big| \L_{X_\al^{(2)}}\L_{X_\al^{(2)}}K \big| \be \big\rangle 
+ \big\langle \al \big|\L_{X_\be^{(2)}}\L_{X_\be^{(2)}}K \big| \al \big\rangle  \\
& \qquad\qquad  -2\big\langle \al \big| \L_{X_\al^{(2)}}\L_{X_\be^{(2)}}K \big| \be \big\rangle \Big)\\
R_{12} &= \int_F \Big( \langle \mathcal D(\al,\al,F),\mathcal F(\be,\be,F)\rangle 
+ \langle \mathcal D(\be,\be,F), \mathcal F(\al,\al,F)\rangle  \\
& \qquad\qquad  -\langle \mathcal D(\al,\be,F)+\mathcal D(\be,\al,F), \mathcal F(\al,\be,F)\rangle \Big) \\
R_2 &= \|\mathcal F(\al,\be,F)\|^2_{K_F}
-\big\langle \mathcal F(\al,\al,F)),\mathcal F(\be,\be,F)\big\rangle_{K_F} \\
R_3 &= -\tfrac34 \| \mathcal D(\al,\be,F)-\mathcal D(\be,\al,F) \|^2_{L_F}
\end{aligned}
}$$
In the case of landmark points, where $m=0, N = \R^D$ and $K$ is diagonal, it is easy to check that our force and stress and the above formula for curvature are exactly the same as those given in our earlier paper \cite{MMM1}. In that paper the individual terms are studied in special cases giving some intuition for them.

\section{\label{nmb:6} Appendix on Convenient Calculus -- Calculus beyond Banach spaces} 
The traditional differential calculus works well for finite dimensional vector spaces and for Banach spaces. 
For more general locally convex spaces we sketch here the convenient approach as explained in \cite{FK88} and \cite{KM97}. 
The main difficulty is that composition of linear mappings stops being jointly continuous at the level of Banach spaces, 
for any compatible topology. 
We use the notation of \cite{KM97} and this is the main reference for the whole appendix.

\subsection{\label{nmb:6.1}Convenient vector spaces and the $c^\infty$-topology}
Let $E$ be a locally convex vector space. A curve $c:\mathbb R\to E$ is called {\it smooth} or $C^\infty$ if all derivatives exist and are 
continuous - this is a concept without problems. Let 
$C^\infty(\mathbb R,E)$ be the space of smooth functions. It can be 
shown that $C^\infty(\mathbb R,E)$ does not depend on the locally convex 
topology of $E$, but only on its associated bornology (system of bounded 
sets).

$E$ is said to be a {\it convenient 
vector space} if one of the following equivalent
conditions is satisfied (called $c^\infty$-completeness):
\begin{enumerate}
\item For any $c\in C^\infty(\mathbb R,E)$ the (Riemann-) integral 
       $\int_0^1c(t)dt$ exists in $E$.
\item A curve $c:\mathbb R\to E$ is smooth if and only if $\la\o c$ is 
       smooth for all $\la\in E'$, where $E'$ is the dual consisting 
       of all continuous linear functionals on $E$.
\item Any Mackey-Cauchy-sequence (i.\ e.\  $t_{nm}(x_n-x_m)\to 0$  
       for some $t_{nm}\to \infty$ in $\mathbb R$) converges in $E$. 
       This is visibly a weak completeness requirement.
\end{enumerate}
The final topology with respect to all smooth curves is called the 
$c^\infty$-topology on $E$, which then is denoted by $c^\infty E$. 
For Fr\'echet spaces it coincides with 
the given locally convex topology, but on the space $\mathcal D$ of test 
functions with compact support on $\mathbb R$ it is strictly finer.

\subsection{\label{nmb:6.2}Smooth mappings} 
Let $E$, $F$, and $G$ be convenient vector spaces, and let $U\subset E$ be $c^\infty$-open. 
Here is the key definition that makes everything work: 
a mapping $f:U\to F$ is called {\it smooth} or $C^\infty$, 
if $f\o c\in C^\infty(\mathbb R,F)$ for all $c\in C^\infty(\mathbb R,U)$.

{\it The main properties of smooth calculus are the following.
\begin{enumerate}
\item For mappings on Fr\'echet spaces this notion of smoothness 
coincides with all other reasonable definitions. Even on 
$\mathbb R^2$ this is non-trivial.
\item Multilinear mappings are smooth if and only if they are 
bounded.
\item If $f:E\supseteq U\to F$ is smooth then the derivative 
$df:U\x E\to F$ is  
smooth, and also $df:U\to L(E,F)$ is smooth where $L(E,F)$ 
denotes the space of all bounded linear mappings with the 
topology of uniform convergence on bounded subsets.
\item The chain rule holds.
\item The space $C^\infty(U,F)$ is again a convenient vector space 
where the structure is given by the obvious injection
$$
C^\infty(U,F) \East{C^\infty(c,\ell)}{} 
\negthickspace\negthickspace\negthickspace\negthickspace\negthickspace
\negthickspace\negthickspace\negthickspace
\prod_{c\in C^\infty(\mathbb R,U), \ell\in F^*} 
\negthickspace\negthickspace\negthickspace\negthickspace\negthickspace
\negthickspace\negthickspace\negthickspace
C^\infty(\mathbb R,\mathbb R),
\quad f\mapsto (\ell\o f\o c)_{c,\ell},
$$
where $C^\infty(\mathbb R,\mathbb R)$ carries the topology of compact 
convergence in each derivative separately.
\item The exponential law holds: For $c^\infty$-open $V\subset F$, 
$$
C^\infty(U,C^\infty(V,G)) \cong C^\infty(U\x V, G)
$$
is a linear diffeomorphism of convenient vector spaces. Note 
that this is the main assumption of variational calculus where a smooth curve in a space of functions is assumed to be just a smooth function in one variable more..
\item A linear mapping $f:E\to C^\infty(V,G)$ is smooth (bounded) if 
and only if $E \East{f}{} C^\infty(V,G) \East{\on{ev}_v}{} G$ is smooth 
for each $v\in V$. This is called the smooth uniform 
boundedness theorem \cite[5.26]{KM97}.
\item The following canonical mappings are smooth.
\begin{align*}
&\operatorname{ev}: C^\infty(E,F)\x E\to F,\quad 
\operatorname{ev}(f,x) = f(x)\\
&\operatorname{ins}: E\to C^\infty(F,E\x F),\quad
\operatorname{ins}(x)(y) = (x,y)\\
&(\quad)^\wedge :C^\infty(E,C^\infty(F,G))\to C^\infty(E\x F,G)\\
&(\quad)^\vee :C^\infty(E\x F,G)\to C^\infty(E,C^\infty(F,G))\\
&\operatorname{comp}:C^\infty(F,G)\x C^\infty(E,F)\to C^\infty(E,G)\\
&C^\infty(\quad,\quad):C^\infty(F,F_1)\x C^\infty(E_1,E)\to 
C^\infty(C^\infty(E,F),C^\infty(E_1,F_1))\\
&\qquad (f,g)\mapsto(h\mapsto f\o h\o g)\\
&\prod:\prod C^\infty(E_i,F_i)\to C^\infty(\prod E_i,\prod F_i)
\end{align*}
\end{enumerate}
}

Smooth mappings are always continuous for the $c^\infty$-topology but there are smooth mappings which are not continuous in the given topology of $E$. This is unavoidable and not so horrible as it might appear at first sight. For example the evaluation $E\x E^*\to\mathbb R$ is jointly continuous if and only if $E$ is normable, but it is always smooth.

\end{document}